\documentclass[11pt]{article}
\usepackage{amssymb,latexsym}
\usepackage{amsmath}

\title{ \large {\bf One proof of the original Kemer's theorems \\
(concerning the text of C. Procesi "What happened to PI-theory" \\ arxiv.org/abs/1403.5673). }}

\author{  \bigskip {\normalsize  I.Sviridova} \\
Departamento de Matem\'atica,\\
Universidade de Bras\'\i lia,\\
70910-900 Bras\'\i lia, DF, Brazil;\\
I.Sviridova@mat.unb.br;\\ sviridova\_i@rambler.ru }

\date{March 07, 2015.}

\newtheorem{theorem}{Theorem}
\newtheorem{definition}{Definition}
\newtheorem{lemma}{Lemma}

\newtheorem{corollary}[lemma]{Corollary}
\newtheorem{remark}{Remark}

\begin{document}
\maketitle

\begin{abstract}
We consider associative algebras over a field of characteristic zero.
We give a version of the proof of the
Kemer's theorems concerning the Specht problem solution \cite{Kem1}, \cite{Kem2}-\cite{Kem5}.
It is proved that the ideal of graded identities
of a finitely generated PI-superalgebra coincides with the
ideal of graded identities of some finite dimensional superalgebra.
This implies that the ideal of
polynomial identities of any (not necessary finitely generated)
PI-algebra coincides with the ideal of identities of the Grassmann envelope
of a finite dimensional superalgebra, and is finitely generated as a
T-ideal.
\bigskip


\end{abstract}

%a partial case of the results of the article published in Comm. Algebra, $G=\mathbb{Z}/2 \mathbb{Z}$

\section*{Introduction}

We present here one of the proofs of the Kemer's theorems about Specht problem solution
for associative PI-algebras over a field of characteristic zero \cite{Kem1}.
We give principally the proof of the fact that any finitely generated $\mathbb{Z}/2 \mathbb{Z}$-graded PI-algebra has the same graded identities as some finite dimensional
$\mathbb{Z}/2 \mathbb{Z}$-graded algebra. This result is the crucial and the most difficult step
of the whole Kemer's solution of the Specht problem.
The proof represented here is the partial case of the proof of a similar result for algebras graded by a finite abelian group given by the author in \cite{Svi1} (the case $G=\mathbb{Z}/2 \mathbb{Z}$).
It is only slightly modified to extend the result from an algebraically closed base field
(as it was considered in \cite{Svi1}) to any field of characteristic zero.
Observe that the proof of generalised results in \cite{Svi1} is independent on the
original Kemer's proof. The unique essential reference to the Kemer's results in that paper
served to obtain Lemma 1. This reference can be safely exchanged by the Lewin's results  \cite{Lewin}.
PI-representability of a non-graded finitely generated PI-algebra can be also obtained as a partial case
of \cite{Svi1} considering the trivial group $G=\{ e \}.$

The second step of the Kemer's arguments is to prove that any PI-algebra
satisfies the same graded polynomial identities as
the Grassmann envelope of some finitely generated $\mathbb{Z}/2 \mathbb{Z}$-graded PI-algebra.
We refer the reader to the book \cite{GZbook} for the proof of this fact.
The author also can prove this fact independently but assume that the proof
represented in \cite{GZbook} is the most elegant, short and clear to understand for a reader.
Thus the author thinks that there is no any need to repeat these arguments.
These two theorems imply (see also \cite {GZbook}) that any PI-algebra
satisfies the same polynomial (non-graded) identities as the Grassmann envelope of some
finite dimensional $\mathbb{Z}/2 \mathbb{Z}$-graded algebra.
The positive solution of the Specht problem in the classical case follows from
the last classification theorem immediately.

Observe that besides the original proof of Aleksandr Kemer \cite{Kem2}-\cite{Kem5},
published later also in \cite{Kem1}, there are several later versions of various authors
including interpretations and generalisations.
Therefore the author does not pretend on originality or some exclusive properties of this text.
And she never would get an idea to rewrite the proof of the original Kemer's results
whenever some obvious misunderstandings appeared some time ago.
This text is just one of the possible version of the proof as the author understood it
generalising these results in \cite{Svi1}.
Notice that this text is a little more detailed then the text published in \cite{Svi1}.
Even here we omit some details and proofs that seems to us not very difficult to restore
for a reader. But the author can assure possible readers that all the statements with omitted
arguments were really checked. Moreover, even the text \cite{Svi1} is the restricted
version of the original detailed author's text.

Throughout the paper we consider only associative algebras
over a field of characteristic zero (not necessary unitary).
Further they will be called algebras.

Let $F$ be a field of characteristic zero
and $F\langle X\rangle$ the free associative algebra over
$F$ generated by a countable set $X=\{x_1,x_2,\ldots \}$.
A T-ideal of $F\langle X\rangle$ is a bilateral ideal invariant
under all endomorphisms of $F\langle X\rangle.$

Let $A$ be an associative algebra over $F$. A polynomial $f=f(x_1,
\ldots , x_n) \in F\langle X\rangle$ is called polynomial
identity for $A$ if $f(a_1, \ldots , a_n)=0$ for any $a_1, \ldots
, a_n \in A$ ($ f \equiv 0$ in $A$). Let us denote by $\mathrm{Id}(A)=\{  f \in
F\langle X\rangle \, | \, f \equiv 0 \, \, \mathrm{in} \, \, A \}$
the ideal of all polynomial identities of $A$. If $A$ satisfies
a non-trivial polynomial identity then $A$ is called $PI$-
algebra. It is well known, for example, that any finite
dimensional algebra is PI (\cite{Jacobs, Proc, Row}). The
relation between T-ideals of $F\langle X\rangle$ and
$PI$-algebras is well understood: for any $F$-algebra $A$,
$\mathrm{Id}(A)$ is a T-ideal of $F\langle X\rangle,$ and any T-
ideal $I$ of $F\langle X\rangle$ is the ideal of identities of
some $F$-algebra $A,$ in particular, of the relatively free algebra
of $F\langle X\rangle/I.$

An algebra $A$ is called $\mathbb{Z}/2 \mathbb{Z}$-graded (or super-algebra) if
$A =A_{\bar{0}} \oplus A_{\bar{1}}$ is the direct sum of its subspaces $A_{\bar{0}}, A_{\bar{1}}$ where
$A_{\theta} A_{\xi} \subseteq A_{\theta \xi}$ holds for any
$\theta, \xi \in \mathbb{Z}/2 \mathbb{Z}.$ An element $a \in A_{\theta}$ is called
homogeneous in the $\mathbb{Z}/2 \mathbb{Z}$-grading of the graded degree $\theta,$ we
also write $\deg_{\mathbb{Z}_2} a =\theta$ in this case. The homogeneous
component $A_{\bar{0}}$ of the $\mathbb{Z}/2 \mathbb{Z}$-graded algebra $A$ is called neutral (or even).

A homomorphism of $\mathbb{Z}/2 \mathbb{Z}$-graded algebras $A,$ $B$ \
$\varphi:A \rightarrow B$ is graded if $\varphi(A_{\theta}) \subseteq B_{\theta}$ for any
$\theta \in \mathbb{Z}/2 \mathbb{Z}.$
An ideal $I \unlhd A$ of a graded algebra $A$
is graded if and only if $I$ is generated by homogeneous in the
grading elements. In this case the quotient algebra $A/I$ is also
$\mathbb{Z}/2 \mathbb{Z}$-graded with the grading induced by the grading of $A$
($\deg_{\theta} \bar{a} = \deg_{\theta} a$).

Let us denote by $A_1 \times
\dots \times A_{\rho}$ the direct product of algebras $A_1,
\dots, A_{\rho},$ and by $A_1 \oplus \dots \oplus A_{\rho}
\subseteq A$ the direct sum of subspaces $A_i$ of an algebra $A.$
We also denote by $J(A)$ the Jacobson radical of
$A.$ Observe that in general all bases and
dimensions of spaces and algebras are defined over the
base field $F$ unless otherwise indicated.

We always assume that the set $\mathbb{N}_0^{k}$
($\mathbb{N}_0=\mathbb{N} \bigcup \{0\}$) is linearly ordered
with the lexicographical order for any natural $k.$

The most part of notions, definitions, facts and
properties of polynomial identities may be found in
\cite{BelRow}, \cite{Dren}, \cite{DrenForm}, \cite{GZbook}, \cite{Kem1}.

\section{Free graded algebra.}

Let us denote by $X^{\mathbb{Z}_2} = \{ x_{i \theta} | i \in \mathbb{N},
\theta \in \mathbb{Z}/2 \mathbb{Z} \}$ a countable set of pairwise different elements.
The algebra $\mathfrak{F}=F\langle X^{\mathbb{Z}_2} \rangle$ is the free associative
$\mathbb{Z}/2 \mathbb{Z}$-graded algebra with the grading $\mathfrak{F}=\oplus_{\theta
\in \mathbb{Z}/2 \mathbb{Z}} \mathfrak{F}_{\theta},$ where
$\mathfrak{F}_{\theta}=\langle x_{i_1 \theta_1}x_{i_2 \theta_2}
\cdots x_{i_s \theta_s} | \theta=\theta_1 + \theta_2 +\cdots +\theta_s
\rangle_{F}.$ Since all the variables are pairwise distinct then
any multihomogeneous polynomial is also a homogeneous element of $\mathfrak{F}$ in
the sense of the grading. An element $x_{i \theta}$ is called graded variables,
and $f \in \mathfrak{F}$ a graded polynomial.

Let $f=f(x_{1 \theta_1}, \dots x_{n \theta_n}) \in F\langle X^{\mathbb{Z}_2}
\rangle$ be a non-trivial $\mathbb{Z}/2 \mathbb{Z}$-graded polynomial. We say that a
$\mathbb{Z}/2 \mathbb{Z}$-graded algebra $A$ satisfies the graded identity $f(x_{1
\theta_1}, \dots x_{n \theta_n})=0$ if $f(a_{1 \theta_1}, \dots
a_{n \theta_n})=0$ in $A$ for any $a_{i \theta_i} \in
A_{\theta_i}.$ Denote by $\mathrm{Id}^{\mathbb{Z}_2}(A) \unlhd
F\langle X^{\mathbb{Z}_2} \rangle$  the ideal of
$\mathbb{Z}/2 \mathbb{Z}$-graded identities of $A.$ Similar
to the case of ordinary (non-graded) identities $\mathrm{Id}^{\mathbb{Z}_2}(A)$
is a two-side graded ideal of  $F\langle X^{\mathbb{Z}_2} \rangle$
invariant under all graded endomorphisms
of $F\langle X^{\mathbb{Z}_2} \rangle$. Such ideals are called
$\mathbb{Z}_2$T-ideals. It is clear that any $\mathbb{Z}_2$T-ideal $I$ of $F\langle
X^{\mathbb{Z}_2}\rangle$ is the ideal of $\mathbb{Z}/2 \mathbb{Z}$-graded identities of
$F\langle X^{\mathbb{Z}_2}\rangle/I.$

Take a set $S \subseteq F\langle X^{\mathbb{Z}_2} \rangle.$
Denote by $\mathbb{Z}_2 T[S]$
the $\mathbb{Z}_2$T-ideal generated by $S.$ Then $\mathbb{Z}_2 T[S]$ contains
exactly all graded identities which are consequences of polynomials
of the set $S.$ $\mathbb{Z}/2 \mathbb{Z}$-graded algebras $A$ and $B$
are called $\mathbb{Z}_2$PI-equivalent ($A \sim_{\mathbb{Z}_2} B$) if
$\mathrm{Id}^{\mathbb{Z}_2}(A)=\mathrm{Id}^{\mathbb{Z}_2}(B).$
Let $\Gamma)$ be a $\mathbb{Z}_2$T-ideal.
We write also $f=g \ (\mathrm{mod } \ \Gamma)$ for $f, g \in F\langle X^{\mathbb{Z}_2} \rangle$ iff $f-g \in \Gamma.$

We consider only graded identities of PI-superalgebras, i.e. the
case when $\mathrm{Id}^{\mathbb{Z}_2}(A) \supseteq \Gamma$ for some nonzero ordinary
T-ideal $\Gamma \unlhd F\langle X \rangle.$ This holds iff the neutral component
$A_{e}$ is a PI-algebra ($\mathrm{Id}^{\mathbb{Z}_2}(A) \ni f(x_{1
\bar{0}},\dots,x_{n \bar{0}}) \ne 0$) (see \cite{BergC}, and \cite{BahtGR}).

Observe that a finitely generated superalgebra always has a
finite set of $\mathbb{Z}/2 \mathbb{Z}$-homogeneous generators.
We consider only homogeneous generating sets of superalgebras.

\begin{definition}
Given a finitely generated superalgebra
$A$ and a finite homogeneous generating set $K$ of $A$
denote by $\mathrm{grk}(K)$ the maximal number of elements
of the same graded degree in $K.$

Then the homogeneous rank $\mathrm{grk}(A)$ of $A$ is the least $\mathrm{grk}(K)$
for all finite homogeneous generating sets $K$ of $A.$
\end{definition}

Let $X_\nu^{\mathbb{Z}_2} = \{ x_{i \theta} | 1 \le i \le \nu; \  \theta \in
\mathbb{Z}/2 \mathbb{Z} \}$ be a finite set of graded variables and
$F\langle X_\nu^{\mathbb{Z}_2} \rangle$ the free associative $\mathbb{Z}/2 \mathbb{Z}$-graded
algebra of the rank $\nu,$ \ $\nu \in \mathbb{N}.$

Let $A$ be a finitely generated $\mathbb{Z}/2 \mathbb{Z}$-graded algebra.
Then  the superalgebra \\$U_\nu=F\langle X_\nu^{\mathbb{Z}_2}
\rangle/(\mathrm{Id}^{\mathbb{Z}_2}(A) \bigcap F\langle X_\nu^{\mathbb{Z}_2} \rangle)$ is the relatively free algebra of the rank $\nu$ for $A.$
Similarly to case of ordinary identities if $\nu \ge
\mathrm{grk}(A)$ then
$\mathrm{Id}^{\mathbb{Z}_2}(A)=\mathrm{Id}^{\mathbb{Z}_2}(U_\nu).$
Moreover the next remark takes place.

Given a $\mathbb{Z}_2$T-ideal $\Gamma_1 \subseteq F\langle X^{\mathbb{Z}_2} \rangle$ denote by
$\Gamma_1(F\langle X_\nu^{\mathbb{Z}_2} \rangle)=\{ f(h_1,\dots,h_n) | f \in
\Gamma_1, h_i \in  F\langle X_\nu^{\mathbb{Z}_2} \rangle,
 \  \deg_{\mathbb{Z}_2} h_i = \deg_{\mathbb{Z}_2} x_i, \forall i \} \unlhd F\langle X_\nu^{\mathbb{Z}_2}
\rangle$ the verbal ideal of the free associative $\mathbb{Z}/2 \mathbb{Z}$-graded
algebra of the rank $\nu$ generated by all appropriate
evaluations of elements of $\Gamma_1.$

\begin{remark} \label{freefg}
\begin{enumerate}
\item $f(x_1,\dots,x_n) \in \mathrm{Id}^{\mathbb{Z}_2}(F\langle
X_\nu^{\mathbb{Z}_2} \rangle/(\Gamma \bigcap F\langle X_\nu^{\mathbb{Z}_2} \rangle))$ if
and only if \\ $f(h_1,\dots,h_n) \in \Gamma$ for all homogeneous
polynomials of appropriate graded degrees $h_1,\dots,h_n \in
F\langle X_\nu^{\mathbb{Z}_2} \rangle.$
\item Let $A$ be a finitely generated associative $\mathbb{Z}/2 \mathbb{Z}$-graded
algebra, $\Gamma_2=\mathrm{Id}^{\mathbb{Z}_2}(A)$ the ideal of its
graded identities, and $\nu \ge \mathrm{grk}(A).$ Then
$\Gamma_1(F\langle X_\nu^{\mathbb{Z}_2}
\rangle) \subseteq \Gamma_2$ implies that $\Gamma_1 \subseteq \Gamma_2.$
\end{enumerate}
\end{remark}

It's well known due to the linearization process and
possibility to identify variables in case of zero
characteristic that any system of identities (ordinary or $\mathbb{Z}/2 \mathbb{Z}$-graded)
is equivalent to a system of multilinear identities. Thus in
case of zero characteristic it is enough to consider only
multilinear identities.

\section{Graded algebras.}

Let $\widetilde{F}$ be an algebraically closed field. Consider a finite dimensional $\widetilde{F}$-superalgebra $A=A_{\bar{0}} \oplus A_{\bar{1}}.$ Then we have an analog of
the Wedderburn-Maltsev decomposition for superalgebras.

\begin{lemma} \label{Pierce}
Let $\widetilde{F}$ be an algebraically closed field.
Any finite dimensional $\widetilde{F}$-super\-al\-geb\-ra $A$ is isomorphic as a
superalgebra to an $F$-superalgebra of the form
\begin{eqnarray} \label{matrix}
A'=(C_1 \times \dots \times
C_p) \oplus J.
\end{eqnarray}
Where the Jacobson radical $J=J(A')$ of $A'$ is a graded ideal,
\  $B'=C_1 \times \dots \times
C_p$ is a maximal graded semisimple
subalgebra of $A',$ \ $p \in \mathbb{N} \bigcup \{0\}.$ A
$\mathbb{Z}/2 \mathbb{Z}$-graded simple component $C_l$ is one of the superalgebras
$M_{k_l,m_l}(\widetilde{F})$ (a matrix algebra with an elementary $\mathbb{Z}/2 \mathbb{Z}$-grading)
or $M_{k_l}(\widetilde{F}[c])$ (a matrix algebra over the group algebra of $\widetilde{F}[c|c^2=1]$ with the grading induced by the natural grading of the group algebra).

Moreover $A'=\mathrm{Span}_{\widetilde{F}}\{ D, U \},$ \ where
\begin{eqnarray} \label{basis}
\label{basisD} &&D=\cup_{1 \le l \le p} \  D_l, \nonumber \\
&&D_l=\{ E_{l i_l j_l} = \varepsilon_l E_{l i_l j_l} \varepsilon_l \
| \  \ 1 \le i_l, j_l \le s_l \
\}, \nonumber \\
&&s_l=k_l+m_l \quad \ \mbox{ if } \ \ C_l=M_{k_l,m_l}(\widetilde{F}); \nonumber \\
&&D_l=\{ E_{l i_l j_l}=\varepsilon_l E_{l i_l j_l} \varepsilon_l, \
E_{l i_l j_l} c =\varepsilon_l \  E_{l i_l j_l} c \  \varepsilon_l \
| \ 1 \le i_l, j_l \le s_l \}, \\
&&s_l=k_l, \  c^2=1 \quad \  \mbox{ if } \  \  C_l=M_{k_l}(\widetilde{F}[c]); \nonumber \\
\label{basisU} &&U=\{ \varepsilon_{l'} u_{\theta}
\varepsilon_{l''} | \ l', l''=1,\dots,p+1; \ u_{\theta} \in
J_{\theta}=J \cap A'_{\theta}, \ \theta \in \mathbb{Z}/2 \mathbb{Z}  \}
\end{eqnarray}
are the sets of homogeneous in the grading elements. $E_{l i_l j_l}$
is the matrix unit. Elements of $D$ are homogeneous with the grading defined by the next equalities
$\deg_{\mathbb{Z}_2} E_{l i_l j_l} = \bar{0}$ if $1 \le i_l, j_l \le k_l,$  or
$k_l+1 \le i_l, j_l \le k_l+m_l,$ and $\deg_{\mathbb{Z}_2} E_{l i_l j_l} = \bar{1}$ otherwise
in $C_l=M_{k_l,m_l}(\widetilde{F});$ \
$\deg_{\mathbb{Z}_2} E_{l i_l j_l} = \bar{0},$ \
$\deg_{\mathbb{Z}_2} E_{l i_l j_l} \cdot c = \bar{1}$  for all $1 \le i_l, j_l \le k_l$
in $C_l=M_{k_l}(\widetilde{F}[c]).$
The element $\varepsilon_l $ is a minimal orthogonal central idempotent of $B',$
homogeneous in the grading, that corresponds to
the unit element of $C_l$ ($l=1,\dots,p$).
Elements $u_{\theta}$ runs on the set of homogeneous
basic elements of the Jacobson radical
$J=J(A)=\oplus_{l',l''=1}^{p+1} \varepsilon_{l'} J \varepsilon_
{l''},$ where $\varepsilon_{p+1}$ is the adjoint idempotent.
\end{lemma}

\noindent {\bf Proof.} It is a classical result (see, e.g., \cite{GZbook} or \cite{BahtZaic2}, \cite{BahtZaicSeg} for a more general case) that a finite dimensional superalgebra $A$ can be represented as $B \oplus J,$ where $J$ is the Jacobson radical of
$A$, that is a graded ideal, and $B$ is a maximal $\mathbb{Z}/2 \mathbb{Z}$-graded semisimple
subalgebra of $A.$ It is also well known that $B \cong \times_{l=1}^{p}
C_l,$ where $C_l$ are simple superalgebras. $C_l$ is isomorphic
either to a matrix superalgebra $M_{k_l,m_l}(F)=M_{k_l+m_l}(F)$ with an elementary grading
\begin{equation*}
\
\begin{array}{rl}
& \quad \  {\scriptstyle k_l} \ \ {\scriptstyle m_l} \qquad  \\
(C_l)_{\bar{0}}=\left\{  \begin{array}{r} {\scriptstyle k_l} \\ {\scriptstyle m_l} \end{array}
\right. &
\left. \left(\begin{tabular}{ c|c }
                    $\ast$ & 0
                    \\   \hline
                    0 & $\ast$ \\
                    \end{tabular}  \right) \right\},
\end{array}
\qquad \
\begin{array}{rl}
& \quad \  {\scriptstyle k_l} \ \ {\scriptstyle m_l} \qquad  \\
(C_l)_{\bar{1}}=\left\{  \begin{array}{r} {\scriptstyle k_l} \\ {\scriptstyle m_l} \end{array}
\right. &
\left. \left(\begin{tabular}{ c|c }
                    0 & $\ast$
                    \\   \hline
                    $\ast$ & 0 \\
                    \end{tabular}  \right) \right\},
\end{array}
\end{equation*}
or to a matrix algebra over the group algebra $M_{k_l}(\widetilde{F}[\mathbb{Z}/2 \mathbb{Z}]) \simeq
M_{k_l}(\widetilde{F}[c|c^2=1])$  with the grading induced by the natural $\mathbb{Z}/2 \mathbb{Z}$-grading
of the group algebra $(C_l)_{\bar{0}}=M_{k_l}(\widetilde{F}),$ \ $(C_l)_{\bar{1}}=M_{k_l}(\widetilde{F}) \cdot c.$

Thus $A'=B' \oplus J$ can be assumed our superalgebra. Any $\mathbb{Z}/2 \mathbb{Z}$-simple component $C_l$
of the semisimple graded subalgebra $B'$
has the unit element of the even degree. It gives the minimal orthogonal
central idempotent $\varepsilon_l$ of the algebra $B'.$
$D$ is a homogeneous basis of $B'.$
It is clear that any element $r \in J$ can be uniquely represented
as a sum of elements of graded subspaces $\varepsilon_{l'}
J \varepsilon_{l''},$ \quad $l',l'' =1,\dots, p+1.$ Moreover
$\varepsilon_l a=0$ (and $a \varepsilon_l=0$) for any $a \in
\varepsilon_{l'} J \varepsilon_{l''}$ with $l \ne l'$ ($l \ne l''$),
$l=1,\dots,p.$  \hfill $\Box$

The next construction is useful to extend the previous result
for any field of zero characteristic in some sense.

Take a superalgebra $B$ (not necessarily without unit).  We
denote by $B^{\#}=B \oplus F \cdot 1_F$ the superalgebra
with the adjoint unit $1_F$.

Let us take a finite dimensional superalgebra
$A=B \oplus J(A)$ with a maximal $\mathbb{Z}/2 \mathbb{Z}$-graded semisimple subalgebra $B$
and the Jacobson radical $J(A).$ Consider a $\mathbb{Z}/2 \mathbb{Z}$-graded subalgebra
$\widetilde{B} \subseteq B,$ and a positive integer number $q.$
Consider the free product $\widetilde{B}^{\#} *_F F\langle X_q^{\mathbb{Z}_2}
\rangle^{\#},$ and define on it the $\mathbb{Z}/2 \mathbb{Z}$-grading
by the equalities $\deg_{\mathbb{Z}_2} (u_1 \cdots u_s)=(\deg_{\mathbb{Z}_2} u_1)
+\dots+(\deg_{\mathbb{Z}_2} u_s),$ where $u_i \in
\widetilde{B}^{\#} \bigcup F\langle X_q^{\mathbb{Z}_2} \rangle^{\#}$ are homogeneous elements. Let
$\widetilde{B}(X_q^{\mathbb{Z}_2})$ be the graded subalgebra of
$\widetilde{B}^{\#} *_F F\langle X_q^{\mathbb{Z}_2} \rangle^{\#}$ generated by the
set $\widetilde{B} \bigcup F\langle X_q^{\mathbb{Z}_2} \rangle.$
Denote by $(X_q^{\mathbb{Z}_2})$
the two-sided graded ideal of $\widetilde{B}(X_q^{\mathbb{Z}_2})$ generated
by the set of variables $X_q^{\mathbb{Z}_2}.$
Particularly,
it is clear that $\widetilde{B}(X_q^{\mathbb{Z}_2})=\widetilde{B} \oplus (X_q^{\mathbb{Z}_2}).$

Given a $\mathbb{Z}_2$T-ideal $\Gamma$ denote by
$\Gamma(\widetilde{B}(X_q^{\mathbb{Z}_2}))$ the two-sided graded verbal
ideal of $\widetilde{B}(X_q^{\mathbb{Z}_2})$ generated by results of all
appropriate evaluations of polynomials from $\Gamma.$
Take any $s \in \mathbb{N},$ and consider the quotient algebra
\begin{eqnarray} \label{FRad}
\mathcal{R}_{q,s}(\widetilde{B},\Gamma)=\widetilde{B}(X_q^{\mathbb{Z}_2})/
(\Gamma(\widetilde{B}(X_q^{\mathbb{Z}_2}))+(X_q^{\mathbb{Z}_2})^s).
\end{eqnarray}
Denote also
$\mathcal{R}_{q,s}(A)=\mathcal{R}_{q,s}(B,\mathrm{Id}^{\mathbb{Z}_2}(A))$ for
$\Gamma=\mathrm{Id}^{\mathbb{Z}_2}(A),$ $\widetilde{B}=B.$

\begin{lemma} \label{Aqs}
Take any natural numbers $q, s$ and a $\mathbb{Z}_2$T-ideal $\Gamma$ such that
$\Gamma \subseteq \mathrm{Id}^{\mathbb{Z}_2}(A).$ The algebra
$\mathcal{R}_{q,s}(\widetilde{B},\Gamma)$ is a finite
dimensional superalgebra with the ideal of graded identities
$\mathrm{Id}^{\mathbb{Z}_2}(\mathcal{R}_{q,s}(\widetilde{B},\Gamma)) \supseteq
\Gamma.$ Moreover
$\mathcal{R}_{q,s}(\widetilde{B},\Gamma)=\overline{B} \oplus
J(\mathcal{R}_{q,s}(\widetilde{B},\Gamma)).$ Here $\overline{B}$ is a maximal
semisimple $\mathbb{Z}/2 \mathbb{Z}$-graded subalgebra of $\mathcal{R}_{q,s}(\widetilde{B},\Gamma),$ and $\overline{B} \cong
\widetilde{B}.$ The Jacobson radical of $\mathcal{R}_{q,s}(\widetilde{B},\Gamma)$
is equal to $(X_q^{\mathbb{Z}_2})/
(\Gamma(\widetilde{B}(X_q^{\mathbb{Z}_2}))+(X_q^{\mathbb{Z}_2})^s),$ and is nilpotent of degree
less or equal to $s.$

If $q \ge \mathrm{grk}(J(A))$, and $s \ge \mathrm{nd}(A)$
then $\mathrm{Id}^{\mathbb{Z}_2}(\mathcal{R}_{q,s}(A)) = \mathrm{Id}^{\mathbb{Z}_2}(A).$
\end{lemma}
\noindent {\bf Proof.} It is clear that
$I=\Gamma(\widetilde{B}(X_q^{\mathbb{Z}_2}))+(X_q^{\mathbb{Z}_2})^s \subseteq (X_q^{\mathbb{Z}_2}),$ and
$\widetilde{B} \bigcap I = (0).$ $I$ is a graded ideal of
$\widetilde{B}(X_q^{\mathbb{Z}_2}).$ Hence
$\mathcal{R}_{q,s}(\widetilde{B},\Gamma)=\widetilde{B}(X_q^{\mathbb{Z}_2})/I$ is a
superalgebra. Then for the canonical homomorphism $\psi:\widetilde{B}(X_q^{\mathbb{Z}_2}) \rightarrow \mathcal{R}_{q,s}(\widetilde{B},\Gamma)$ we obtain
$\overline{B}=\psi(\widetilde{B}) \cong \widetilde{B}.$ Hence
$\mathcal{R}_{q,s}(\widetilde{B},\Gamma)=\overline{B} \oplus
\psi((X_q^{\mathbb{Z}_2})),$ where $\psi((X_q^{\mathbb{Z}_2}))=(X_q^{\mathbb{Z}_2})/I$ is the
maximal nilpotent ideal of the algebra
$\mathcal{R}_{q,s}(\widetilde{B},\Gamma)$ of degree at most $s.$
It is clear that $(X_q^{\mathbb{Z}_2})/I$ is $\mathbb{Z}/2 \mathbb{Z}$-graded and finite dimensional.
Then $\mathcal{R}_{q,s}(\widetilde{B},\Gamma)$ is also a finite dimensional
algebra with the Jacobson radical
$J(\mathcal{R}_{q,s}(\widetilde{B},\Gamma))=\psi((X_q^{\mathbb{Z}_2})).$ It is
clear that $\Gamma \subseteq
\mathrm{Id}^{\mathbb{Z}_2}(\mathcal{R}_{q,s}(\widetilde{B},\Gamma))$ for any $q,
s \in \mathbb{N}.$

Let us take $\Gamma=\mathrm{Id}^{\mathbb{Z}_2}(A),$ $\widetilde{B}=B.$
Suppose that the Jacobson radical $J(A)$ of the algebra $A$ is generated
as an algebra by the set $\{r_1,\dots,r_\nu\}$.
Then we have $r_i=\sum_{\theta \in \mathbb{Z}/2 \mathbb{Z}} r_{i \theta},$ where $r_{i \theta} \in
J(A) \bigcap A_{\theta}$ ($i=1,\dots,\nu$).
Consider the map $\varphi:\overline{x}_{i \theta}=x_{i
\theta}+I \mapsto r_{i \theta}$ \  ($i=1,\dots,\nu$).
Assume that $\varphi(b+I)=b$ for any $b \in B.$
If $q \ge \nu,$ and $s \ge
\mathrm{nd}(A)$ then $\varphi$ can be extended to
a surjective graded homomorphism $\varphi:\mathcal{R}_{q,s}(A)
\rightarrow A.$ Therefore
$\Gamma = \mathrm{Id}^{\mathbb{Z}_2}(A) \supseteq
\mathrm{Id}^{\mathbb{Z}_2}(\mathcal{R}_{q,s}(A)).$ \hfill $\Box$

\begin{definition}
$F$-algebra $A$ is called representable if $A$ can be embedded
into some algebra $C$ that is finite dimensional over an extension
$\widetilde{F} \supseteq F$ of the base field $F.$
\end{definition}

\begin{lemma} \label{Repr}
If a $\mathbb{Z}/2 \mathbb{Z}$-graded $F$-algebra $A$ is representable then there exists
a finite dimensional over $F$ $\mathbb{Z}/2 \mathbb{Z}$-graded $F$-algebra $U$ such that
$\mathrm{Id}^{\mathbb{Z}_2}(A)=\mathrm{Id}^{\mathbb{Z}_2}(U).$
\end{lemma}
\noindent {\bf Proof.} Suppose that $A$ is isomorphic to an
$F$-subalgebra $\mathcal{B}$ of a finite dimensional
$\widetilde{F}$-algebra $\widetilde{\mathcal{B}}.$ We can assume
that the extension $\widetilde{F} \supseteq F$ is algebraically
closed, and $\mathcal{B}$ is $\mathbb{Z}/2 \mathbb{Z}$-graded.
Consider the algebra $\widetilde{U}=\widetilde{U}_{\bar{0}} \oplus \widetilde{U}_{\bar{1}},$ where $\widetilde{U}_\theta=(\widetilde{F}
\mathcal{B}_\theta) \otimes_{\widetilde{F}} \widetilde{F} \theta \subseteq
\widetilde{\mathcal{B}} \otimes_{\widetilde{F}} \widetilde{F}[\mathbb{Z}/2 \mathbb{Z}].$
$\widetilde{U}$
is a finite dimensional $\mathbb{Z}/2 \mathbb{Z}$-graded $\widetilde{F}$-algebra. And
$\mathrm{Id}^{\mathbb{Z}_2}(\widetilde{U})=\mathrm{Id}^{\mathbb{Z}_2}(\mathcal{B})=
\mathrm{Id}^{\mathbb{Z}_2}(A).$
Let $\widetilde{C}_l=\widetilde{C}_{l \bar{0}} \oplus \widetilde{C}_{l \bar{1}}$
be a $\mathbb{Z}/2 \mathbb{Z}$-simple component in the decomposition (\ref{matrix}) of
the algebra $\widetilde{U}.$

It is clear from the classification of simple superalgebras (see also Lemma \ref{Pierce})
that $\widetilde{C}_l$ contains a
finite dimensional over $F$ $\mathbb{Z}/2 \mathbb{Z}$-graded simple $F$-subalgebra
$C_l=C_{l \bar{0}} \oplus C_{l \bar{1}}$ satisfying
$\widetilde{C}_{l \theta}=\widetilde{F} C_{l \theta}$ for any
$\theta \in \mathbb{Z}/2 \mathbb{Z}$ ($l=1,\dots,p$).
More precisely, $C_l=M_{k_l,m_l}(F)$ if $\widetilde{C}_l=M_{k_l,m_l}(\widetilde{F}),$
and $C_l=M_{k_l}(F[c])$ if $\widetilde{C}_l=M_{k_l}(\widetilde{F}[c]).$
Moreover the algebras $\widetilde{C}_l$ and $C_l$ have the same canonic homogeneous base
of type (\ref{basisD}) over $\widetilde{F}$ and $F$ respectively.

Let us take $B=C_1 \times \dots \times C_p,$
$\Gamma=\mathrm{Id}^{\mathbb{Z}_2}(A),$ $q=\dim_{\widetilde{F}}
J(\widetilde{U}),$ $s=\mathrm{nd}(\widetilde{U}).$ Then the $F$-algebra
$U=\mathcal{R}_{q,s}(B,\Gamma)$ defined by (\ref{FRad}) is
$\mathbb{Z}/2 \mathbb{Z}$-graded finite dimensional over $F.$ And $\mathrm{Id}^{\mathbb{Z}_2}(U) =
\Gamma= \mathrm{Id}^{\mathbb{Z}_2}(\widetilde{U})$ (Lemma \ref{Aqs}). \hfill
$\Box$

\begin{definition}\label{el-dec}
We say that an $F$-finite dimensional superalgebra $A'$ has an
{\it elementary decomposition } if it satisfies
the assertion of Lemma \ref{Pierce}.
\end{definition}

It is clear that the direct product of superalgebras with elementary decomposition
is the superalgebra with elementary decomposition. Also if $F$ is algebraically closed
then any finite dimensional $F$-superalgebra has an elementary decomposition.

\begin{corollary} \label{eldec-fd}
Let $F$ be a field of characteristic zero.
Any finite dimensional $F$-superalgebra
is $\mathbb{Z}_2$PI-equivalent to a finite dimensional $F$-superalgebra with
elementary decomposition.
\end{corollary}
\noindent {\bf Proof.}
A finite dimensional $F$-superalgebra
$A$ can be naturally embedded to the superalgebra $\widetilde{A}=A \otimes_F
\widetilde{F}$ preserving graded identities. We assume here that
$\deg_{\mathbb{Z}_2} a \otimes \alpha= \deg_G a,$
for all $a \in A,$ \ $\alpha \in \widetilde{F}.$ The superalgebra $\widetilde{A}$
is finite dimensional over $\widetilde{F}.$ By Lemma \ref{Repr} there
exists a finite dimensional $F$-superalgebra $A'$ with elementary decomposition
such that $\mathrm{Id}^{\mathbb{Z}_2}(A') =
\mathrm{Id}^{\mathbb{Z}_2}(\widetilde{A})=\mathrm{Id}^{\mathbb{Z}_2}(A).$
Where all identities are considered
over the field $F.$ \hfill $\Box$

Therefore considering graded identities of a finite dimensional
superalgebra $A$ we always can assume that $A$ has a form
(\ref{matrix}), and a basis of $A$ can be chosen in the set $D
\bigcup U.$ Particularly, for multilinear graded polynomials it is
enough to consider only evaluations by elements of the set $D
\bigcup U.$ Such evaluation of a multilinear graded polynomial
is called {\it elementary}. Elements of the set $D$ are called
semisimple, and elements of the set $U$ are radical.

\begin{definition}
Let $A=B \oplus J$ be a finite dimensional superalgebra,
$B=\oplus_{\theta \in
\mathbb{Z}/2 \mathbb{Z}} B_{\theta}$ a maximal semisimple $\mathbb{Z}/2 \mathbb{Z}$-graded subalgebra of $A,$
and $J(A)=J$ the Jacobson radical of $A.$  We denote by
$\mathrm{dims}_{\mathbb{Z}_2} A=(\dim B_{\bar{0}},\dim
B_{\bar{1}}),$ and by $\mathrm{nd}(A)$ the nilpotency degree of the radical $J.$
Consider as principal the next parameter of $A$  \
$\mathrm{par}_{\mathbb{Z}_2}(A)=(\mathrm{dims}_{\mathbb{Z}_2} A;\mathrm{nd}(A)).$

The $4$-tuple $\mathrm{cpar}_{\mathbb{Z}_2}(A)=(\mathrm{par}_{\mathbb{Z}_2}(A);\dim
J(A))$ is called complex parameter of $A.$
\end{definition}
A finite dimensional superalgebra $A$ is
nilpotent if and only if $\mathrm{dims}_{\mathbb{Z}_2} A=(0,0).$
Recall that $n$-tuples of numbers are ordered lexicographically. Then for any
nonzero graded two-side ideal $I \unlhd A$ we have that
$\mathrm{cpar}_{\mathbb{Z}_2}(A/I) < \mathrm{cpar}_{\mathbb{Z}_2}(A).$

\section{Kemer index.}

Let $f=f(y_1,\dots,y_k,x_1,\dots,x_n) \in F\langle X^{\mathbb{Z}_2} \rangle$
be a polynomial linear in all variables of the set $Y=\{ y_1,
\dots, y_k\}.$ The polynomial $f$ is alternating in $Y,$ if
\[f(y_{\sigma(1)},\dots,y_{\sigma(k) },x_1,\dots,x_n)=(-1)^{\sigma}
f(y_1,\dots,y_k,x_1,\dots,x_n)\] holds for any permutation $\sigma \in
\mathrm{Sym}_k.$

For any polynomial $g(y_1,\dots,y_k,x_1,\dots,x_n)$
that is linear in $Y=\{ y_1, \dots, y_k\},$ it is possible to
construct a polynomial alternating in $Y$ by setting
$$f(y_1,\dots,y_k,x_1,\dots,x_n)=
\mathcal{A}_{Y}(g)=\sum_{\sigma \in \mathrm{Sym}_k} (-1)^{\sigma}
g(y_{\sigma(1)},\dots,y_{\sigma(k)},x_1,\dots,x_n).$$ The
corresponding mapping $\mathcal{A}_{Y}$ is a linear transformation, we called it
the alternator. Any polynomial $f$ alternating in $Y$
can be decomposed as $f=\sum_{i=1}^{s} \alpha_i
\mathcal{A}_{Y}(u_i),$ where  $u_i$s are monomials of $f,$ \
$\alpha_i \in F.$ The properties of graded alternating polynomials
are similar to that of non-graded polynomials.
We consider only polynomials alternating in homogeneous sets of variables.

Given a pair $\overline{t}=(t_1,t_2) \in
\mathbb{N}_0^2$ we say that a graded polynomial $f \in F\langle
X^{\mathbb{Z}_2} \rangle$ has a collection of $\overline{t}$-alternating
graded variables ($f$ is $\overline{t}$-alternating) if
$f(Y_{\bar{0}},Y_{\bar{1}},X)$ is linear in
$Y=Y_{\bar{0}} \cup Y_{\bar{1}},$ and $f$ is alternating in each
set $Y_{\theta} = \{y_{1 \theta},\dots,y_{t_{\theta}
\theta} \} \subseteq X_{\theta},$ \
$|Y_{\theta}|=t_{\theta},$ \ $\theta \in \mathbb{Z}/2 \mathbb{Z}.$

\begin{definition}
Fix any $\overline{t}=(t_{\bar{0}},t_{\bar{1}}) \in
\mathbb{N}_0^2.$ Suppose that $\tau_1, \dots, \tau_s \in \mathbb{N}_0^2$
are some (possibly different) pairs
satisfying the conditions $\tau_j=(\tau_{j \bar{0}},\tau_{j \bar{1}})
> \overline{t}=(t_{\bar{0}},t_{\bar{1}})$ for all $j=1, \dots, s.$
Let $f \in F\langle X^{\mathbb{Z}_2} \rangle$
be a multihomogeneous graded polynomial.
Suppose that $f=f(Y_1,\dots,Y_{s+\mu};X)$ has $s$
collections of $\tau_j$-alternating variables $Y_j=Y_{j \bar{0}} \cup Y_{j \bar{1}}$
($j=1,\dots,s$), and $\mu$ collections of $\overline{t}$-alternating
variables $Y_j=Y_{j \bar{0}} \cup Y_{j \bar{1}}$ \
($j=s+1,\dots,s+\mu$). We assume that all these sets are
disjoint. Then we say that $f$ is of the type $(\overline{t};s;\mu)=(t_1,t_2;s;\mu).$
Here $Y_{j \theta} \subseteq
X_{\theta}$ with $|Y_{j \theta}|=\tau_{j \theta}$ for
any $j=1,\dots,s,$ and $|Y_{j \theta}|=t_{\theta}$ \  for any
$j=s+1,\dots,s+\mu.$
\end{definition}

Observe that a multihomogeneous polynomial f of a type $(t;s;\mu)$ is also of the type
$(t;s';\mu')$ for all $s' \le s,$ and $\mu' \le \mu.$
Particularly, any nontrivial graded multilinear polynomial of
degree $s$ has the type $(0,0;s;\mu)$ for any $\mu \in
\mathbb{N}_0.$

\begin{definition} \label{defbeta}
Given a $\mathbb{Z}_2$T-ideal $\Gamma \unlhd F\langle X^{\mathbb{Z}_2} \rangle$ the parameter
$\beta(\Gamma)=(t_{\bar{0}},t_{\bar{1}})$ is the greatest lexicographic
pair $\overline{t}=(t_{\bar{0}},t_{\bar{1}}) \in \mathbb{N}_0^2$ such
that for any $s \in \mathbb{N}$ there exists a graded polynomial
$f \notin \Gamma$ of the type $(\overline{t};0;s).$
\end{definition}

The parameter $\beta(\Gamma)$ is well defined for any proper $\mathbb{Z}_2$T-ideal $\Gamma
\unlhd F\langle X^{\mathbb{Z}_2} \rangle$ of  a
finitely generated $\mathbb{Z}/2 \mathbb{Z}$-graded PI-algebra. In this case
$\Gamma$ contains the ordinary Capelli polynomial of some order $d$
(\cite{Kem0}). Hence any graded polynomial $f$ of the type
$(t_{\bar{0}},t_{\bar{1}};0;s)$  belongs to $\Gamma$
if $t_{\bar{0}} \ge d$ or $t_{\bar{1}} \ge d.$ The next parameter is also well defined.

\begin{definition} \label{defgamma}
Given a nonnegative integer $\mu$ let $\gamma(\Gamma;\mu)=s \in
\mathbb{N}$ be the smallest integer $s>0$ such that any graded
polynomial of the type $(\beta(\Gamma);s;\mu)$ belongs to
$\Gamma.$

$\gamma(\Gamma;\mu)$ is a positive non-increasing function of
$\mu.$ Let us denote the limit of this function by $\gamma(\Gamma)=\lim \limits_{\mu \to
\infty} \gamma(\Gamma;\mu) \in \mathbb{N}.$ Then $\omega(\Gamma)$ is
the smallest number $\widehat{\mu}$ such that
$\gamma(\Gamma;\mu)=\gamma(\Gamma)$ for any $\mu \ge
\widehat{\mu}$.
\end{definition}

\begin{definition}
We call by the Kemer index of a $\mathbb{Z}_2$T-ideal $\Gamma$ the
lexicographically ordered collection
$\mathrm{ind}_{\mathbb{Z}_2}(\Gamma)=(\beta(\Gamma); \gamma(\Gamma)).$
\end{definition}
Notice that $\mathrm{ind}_{\mathbb{Z}_2}(\Gamma)$ is greater than $(0,0;1)$ for any proper $\mathbb{Z}_2$T-ideal
$\Gamma.$ We assume also that $\mathrm{ind}_{\mathbb{Z}_2}(F\langle X^{\mathbb{Z}_2} \rangle)=(0,0;1).$

Let us denote
$\omega(A)=\omega(\mathrm{Id}^{\mathbb{Z}_2}(A)),$ \
$\gamma(A;\mu)=\gamma(\mathrm{Id}^{\mathbb{Z}_2}(A);\mu),$ \
$\mathrm{ind}_{\mathbb{Z}_2}(A)=\mathrm{ind}_{\mathbb{Z}_2}(\mathrm{Id}^{\mathbb{Z}_2}(A))$
for a finitely generated PI-superalgebra $A.$

It is clear that $A$ is a nilpotent superalgebra of class
$s$ if and only if $\mathrm{ind}_{\mathbb{Z}_2}(A)=(0,0;s).$ Particularly
$\mathrm{ind}_{\mathbb{Z}_2}(A)=\mathrm{par}_{\mathbb{Z}_2}(A)$ for a nilpotent algebra $A.$
In general case we have
\begin{lemma} \label{bound}
$\mathrm{ind}_{\mathbb{Z}_2}(A) \le \mathrm{par}_{\mathbb{Z}_2} (A)$
for any finite dimensional  superalgebra $A.$
\end{lemma}
\noindent {\bf Proof.} Let us denote $\mathrm{dims}_{\mathbb{Z}_2} A=(t_{\bar{0}},t_{\bar{1}}).$
Suppose that $\beta_{\tilde{\theta}} > t_{\tilde{\theta}}$ for some $\tilde{\theta}=\bar{0}, \bar{1},$ and a
multilinear polynomial $f$ has the type
$(\beta_{\bar{0}},\beta_{\bar{1}};0;\mathrm{nd}(A)).$ Then for any
$j=1,\dots,\mathrm{nd}(A)$ all $\tilde{\theta}$-variables of the alternating set
$Y_{j \tilde{\theta}}$ of $f$ can not be
evaluated only by se\-mi\-simple elements with nonzero result.
Therefore $f \in \mathrm{Id}^{\mathbb{Z}_2}(A),$ and $\mathrm{ind}_{\mathbb{Z}_2}(A) \le
\mathrm{par}_{\mathbb{Z}_2}(A).$ \hfill $\Box$

\begin{definition} \label{Kemerpolyn}
Given a $\mathbb{Z}_2$T-ideal $\Gamma$, and any $\mu \in
\mathbb{N}_0$ a multihomogeneous polynomial $f \in F\langle X^{\mathbb{Z}_2}
\rangle$ is called $\mu$-boundary for $\Gamma$ if $f \notin \Gamma,$ and $f$ has the type
$(\beta(\Gamma);\gamma(\Gamma)-1;\mu).$

Let us denote by $S_\mu(\Gamma)$ the set of all $\mu$-boundary polynomials
for $\Gamma$. Denote also \ $S_\mu(A)=S_\mu(\mathrm{Id}^{\mathbb{Z}_2}(A)),$ \
$K_\mu(\Gamma)=\mathbb{Z}_2 T[S_\mu(\Gamma)],$ \ $K_{\mu,
A}=K_\mu(\mathrm{Id}^{\mathbb{Z}_2}(A))=\mathbb{Z}_2 T[S_\mu(A)].$
\end{definition}

Observe that if the Kemer index is well defined for a $\mathbb{Z}_2$T-ideal $\Gamma$
then $\Gamma$ has multilinear boundary polynomials for all $\mu \in \mathbb{N}_0.$
Moreover a polynomial $f$ belongs to
$S_{\mu}(\Gamma)$ if and only if its full multilinearization
$\widetilde{f}$ belongs to $S_{\mu}(\Gamma).$ Definitions
\ref{defbeta}, \ref{defgamma} immediately imply the next basic properties of the Kemer index
and boundary polynomials.

\begin{lemma} \label{ind1}
Given $\mathbb{Z}_2$T-ideals $\Gamma_1,$ $\Gamma_2$ admitting the Kemer index \  if $\Gamma_1 \subseteq \Gamma_2$ then $\mathrm{ind}_{\mathbb{Z}_2}(\Gamma_1) \ge \mathrm{ind}_{\mathbb{Z}_2}(\Gamma_2).$
\end{lemma}

\begin{lemma} \label{ind2}
Consider $\mathbb{Z}_2$T-ideals $\Gamma,$ \ $\Gamma_1,$ $\dots,$ $\Gamma_{\rho}$ admitting the Kemer index.
Assume that $\mathrm{ind}_{\mathbb{Z}_2}(\Gamma_i) <
\mathrm{ind}_{\mathbb{Z}_2}(\Gamma)$ for all $i=1,\dots,\rho.$ Then there exists  $\widehat{\mu} \in
\mathbb{N}_0$ such that $S_\mu(\Gamma) \subseteq \bigcap
\limits_{i=1} \limits^{\rho} \Gamma_i$ for any $\mu \ge
\widehat{\mu}.$
\end{lemma}
\noindent {\bf Proof.} Let us denote
$\mathrm{ind}_{\mathbb{Z}_2}(\Gamma_i)=(\beta_i,\gamma_i),$ \ $1 \le i \le \rho,$
\ $\mathrm{ind}_{\mathbb{Z}_2}(\Gamma)=(\beta,\gamma).$ Assume that $\beta_i < \beta$
for $i=1, \dots, \rho',$ and $\beta_i=\beta,$ \ $\gamma_i<\gamma$
for $i=\rho'+1, \dots, \rho$ ($0 \le \rho' \le \rho$).

If $\beta > \beta_i$ ($1 \le i \le \rho'$) then there exists $\mu_i$ such
that any polynomial of the type $(\beta;0;\mu_i)$ belongs to
$\Gamma_i.$ If $\beta_i=\beta$ and $\gamma_i<\gamma$
($i=\rho'+1,\dots,\rho$) then any polynomial of the type
$(\beta;\gamma_i;\mu_i)$ belongs to $\Gamma_i$ for all $\mu_i \ge
\omega(\Gamma_i).$ Thus $S_\mu(\Gamma) \subseteq \bigcap
\limits_{i=1} \limits^{\rho} \Gamma_i$ for any $\mu \ge \max \{
\mu_1,\dots,\mu_{\rho'},\omega(\Gamma_{\rho'+1}),\dots,\omega(\Gamma_{\rho})\}.$
\hfill $\Box$

Lemmas \ref{ind1}, \ref{ind2} jointly give the next properties.

\begin{lemma} \label{ind11}
Given $\mathbb{Z}_2$T-ideals $\Gamma_1,$ $\Gamma_2$ admitting the Kemer index \
$\mathrm{ind}_{\mathbb{Z}_2}(\Gamma_1 \cap \Gamma_2) = \max_{i=1,2} \mathrm{ind}_{\mathbb{Z}_2}(\Gamma_{i}).$
\end{lemma}

\begin{lemma} \label{dirsum}
For all finitely generated PI-superalgebras $A_i$ \
$\mathrm{ind}_{\mathbb{Z}_2}(A_1 \times \dots \times A_{\rho})=\max \limits_{1
\le i \le \rho} \mathrm{ind}_{\mathbb{Z}_2}(A_i).$
\end{lemma}

\begin{lemma} \label{subset}
Given $\mathbb{Z}_2$T-ideals $\Gamma_1,$ $\Gamma_2$ admitting the Kemer index, and satisfying
$\Gamma_1 \subseteq \Gamma_2$  one of the following
alternatives takes place:
\begin{enumerate}
  \item $\mathrm{ind}_{\mathbb{Z}_2}(\Gamma_1) = \mathrm{ind}_{\mathbb{Z}_2}(\Gamma_2),$ \
  and \ $S_\mu(\Gamma_1) \supseteq S_\mu(\Gamma_2),$ \ $K_\mu(\Gamma_1) \supseteq
  K_\mu(\Gamma_2)$ \ $\forall \mu \in \mathbb{N}_0;$
  \item $\mathrm{ind}_{\mathbb{Z}_2}(\Gamma_1) > \mathrm{ind}_{\mathbb{Z}_2}(\Gamma_2),$ \
and \ $S_{\hat{\mu}}(\Gamma_1)
  \subseteq \Gamma_2$ for some $\hat{\mu} \in \mathbb{N}_0.$
\end{enumerate}
Moreover in the case $\Gamma_1 \subseteq \Gamma_2$ the conditions
\  $\mathrm{ind}_{\mathbb{Z}_2}(\Gamma_1) > \mathrm{ind}_{\mathbb{Z}_2}(\Gamma_2)$ and
$S_\mu(\Gamma_1) \subseteq \Gamma_2$  are equivalent for some $\mu
\in \mathbb{N}_0.$
\end{lemma}
\noindent {\bf Proof.} By Lemma \ref{ind1}
we have that $\mathrm{ind}_{\mathbb{Z}_2}(\Gamma_1) \ge \mathrm{ind}_{\mathbb{Z}_2}(\Gamma_2).$ The
conditions $\mathrm{ind}_{\mathbb{Z}_2}(\Gamma_1) =
\mathrm{ind}_{\mathbb{Z}_2}(\Gamma_2)=(\beta,\gamma),$ $\Gamma_1 \subseteq
\Gamma_2$ imply that $S_\mu(\Gamma_1) \supseteq S_\mu(\Gamma_2),$ \
$K_\mu(\Gamma_1) \supseteq K_\mu(\Gamma_2)$ \ $\forall \mu \in
\mathbb{N}_0.$

If $\mathrm{ind}_{\mathbb{Z}_2}(\Gamma_1) > \mathrm{ind}_{\mathbb{Z}_2}(\Gamma_2)$ then
$S_{\hat{\mu}}(\Gamma_1) \subseteq \Gamma_2$ holds for
some $\hat{\mu} \in \mathbb{N}_0$ by Lemma \ref{ind2}. And visa versa if
$\Gamma_1 \subseteq \Gamma_2,$ and $S_{\hat{\mu}}(\Gamma_1) \subseteq
\Gamma_2$ for some $\hat{\mu} \in \mathbb{N}_0$ then
$S_{\hat{\mu}}(\Gamma_2) \subseteq S_{\hat{\mu}}(\Gamma_1) \subseteq
\Gamma_2$ gives a contradiction. Therefore in this case we obtain that $\mathrm{ind}_{\mathbb{Z}_2}(\Gamma_1)
> \mathrm{ind}_{\mathbb{Z}_2}(\Gamma_2).$ \hfill $\Box$

Similarly the conditions $\mathrm{ind}_{\mathbb{Z}_2}(\Gamma_1) = \mathrm{ind}_{\mathbb{Z}_2}(\Gamma_2)$
and  $S_\mu(\Gamma_1) \supseteq S_\mu(\Gamma_2)$ are also equivalent
in the case $\Gamma_1 \subseteq \Gamma_2.$

The last lemma has the following corollary.

\begin{lemma} \label{dirsumKP}
Given an integer $\mu
\in \mathbb{N}_0$ and finitely generated PI-superalgebras
$A_1,\dots,A_\rho$ with the same Kemer index \ $\mathrm{ind}_{\mathbb{Z}_2}
(A_i) =\kappa$ for all $i=1,\dots,\rho$ it holds
\begin{eqnarray*}
&&S_\mu(A_1 \times \dots \times A_\rho)=S_\mu(\bigcap \limits_{i=1}
\limits^{\rho} \mathrm{Id}^{\mathbb{Z}_2}(A_i))=\bigcup \limits_{i=1}
\limits^{\rho} S_\mu(A_i), \\ &&K_{\mu, A_1 \times \dots \times
A_\rho}=K_\mu(\bigcap \limits_{i=1} \limits^{\rho}
\mathrm{Id}^{\mathbb{Z}_2}(A_i))=\sum \limits_{i=1} \limits^{\rho}
K_\mu(\mathrm{Id}^{\mathbb{Z}_2}(A_i))=\sum \limits_{i=1} \limits^{\rho} K_{\mu,
A_i}.
\end{eqnarray*}
\end{lemma}

\begin{lemma} \label{addKP}
Given a $\mathbb{Z}_2$T-ideal $\Gamma,$ and a non-positive integer $\mu \in \mathbb{N}_0$ we have that
$\mathrm{ind}_{\mathbb{Z}_2} (\Gamma) > \mathrm{ind}_{\mathbb{Z}_2} (\Gamma+K_\mu(\Gamma)).$
\end{lemma}
\noindent {\bf Proof.} Since $\Gamma \subseteq \Gamma+K_\mu(\Gamma)$
and $S_\mu(\Gamma) \subseteq K_\mu(\Gamma) \subseteq
\Gamma+K_\mu(\Gamma)$ then the assertion immediately follows from Lemma
\ref{subset}. \hfill $\Box$

\section{$\mathbb{Z}_2$PI-reduced algebras.}

\begin{definition}
A finite dimensional superalgebra $A$ with elementary decomposition is
$\mathbb{Z}_2$PI-reduced if there do not exist finite dimensional
superalgebras $A_1,\dots,A_\varrho$ with elementary decomposition such
that $\bigcap\limits_{i=1}\limits^{\varrho} \mathrm{Id}^{\mathbb{Z}_2}(A_i) =
\mathrm{Id}^{\mathbb{Z}_2}(A),$ and $\mathrm{cpar}_{\mathbb{Z}_2}(A_i) < \mathrm{cpar}_{\mathbb{Z}_2}(A)$
for all $i=1,\dots,\varrho.$
\end{definition}
It is clear that a nilpotent finite dimensional superalgebra
$A$ is $\mathbb{Z}_2$PI-reduced if and only if it has the minimal dimension
among all nilpotent finite dimensional superalgebras satisfying the
same graded identities as $A.$

\begin{lemma} \label{ind-simple}
Any simple finite dimensional superalgebra $A$ with elementary decomposition
is $\mathbb{Z}_2$PI-reduced, and
$\mathrm{ind}_{\mathbb{Z}_2}(A)=\mathrm{par}_{\mathbb{Z}_2}(A)=(t_{\bar{0}},t_{\bar{1}};1)$ (i.e.
$\beta(A)=\mathrm{dims}_{\mathbb{Z}_2} A,$ and $\gamma(A)=\mathrm{nd}(A)=1$).
\end{lemma}
\noindent {\bf Proof.} Suppose that $\mathrm{dims}_{\mathbb{Z}_2} A=(t_{\bar{0}},t_{\bar{1}}).$
Any finite dimensional
$\mathbb{Z}/2 \mathbb{Z}$-graded simple algebra is semisimple, $\mathrm{nd}(A)=1,$ and
$\mathrm{cpar}_{\mathbb{Z}_2}(A)=(\mathrm{dims}_{\mathbb{Z}_2} A;1;0).$
It follows from Lemma \ref{bound} that
$\beta(A) \le \mathrm{dims}_{\mathbb{Z}_2} A=(t_{\bar{0}},t_{\bar{1}}).$ Hence it is
enough to construct a polynomial of
the type $(t_{\bar{0}},t_{\bar{1}};0;\hat{s})$ that is not a graded
identity of $A$ for any $\hat{s} \in \mathbb{N}.$
Fix any natural number $\hat{s}.$
Since $A$ has an elementary decomposition then $A=M_{k,m}(F)$ or $A=M_{k}(F[c]).$

Consider the case $A=M_{k,m}(F)$ at first. In this case
$A$  has a homogeneous in the grading basis of the type $\{ E_{ij} | i, j = 1,\dots,k+m \},$
where $E_{ij}$ are the matrix units.
Consider $\hat{s}$ sets of $(k+m)^2$ distinct graded variables $Y_d=\{y_{d,(ij)} \in X^{\mathbb{Z}_2}|
i,j=1,\dots,k+m \}.$ Here any graded variable $y_{d,(ij)}$
corresponds to the matrix unit $E_{ij},$ and
$\deg_{\mathbb{Z}_2} y_{d,(ij)} = \deg_{\mathbb{Z}_2} E_{ij}$ in the superalgebra $M_{k,m}(F)$
($d=1,\dots,\hat{s}$).

Denote by $Y_{d, \theta}=\{y \in Y_d | \deg{\mathbb{Z}_2} y=\theta \}.$
Then for any $d$ we have that
$|Y_{d, \theta}|=t_{\theta}=\dim A_{\theta},$ \ $\theta \in \mathbb{Z}/2 \mathbb{Z}.$
Let us take also the set of graded variables
$Z=\{z_{(ij)} \in X^{\mathbb{Z}_2} | \deg_{\mathbb{Z}_2} z_{(ij)}=\deg_{\mathbb{Z}_2} E_{ij}, \
i,j=1,\dots,k+m \}.$ We assume that $Z$ is disjoint with $Y=\bigcup_{d=1}^{\hat{s}} Y_d.$
We say that the variable $z_{(j_1 i_2)}$ connects the variables
$y_{d,(i_1 j_1)}$ and $y_{d,(i_2 j_2)}.$ Let us
consider for any fixed $d$ the graded monomial $w_{d}$
that is the product of all variables $y_{d,(i_s j_s)}$
connected by variables $z_{(ij)}$
\begin{eqnarray} \label{word1}
&&w_{d}= y_{d,(11)} z_{(11)} y_{d,(12)} z_{(21)}
y_{d,(13)} \cdots  y_{d,(i_1 j_1)} z_{(j_1 i_2)}
y_{d,(i_2 j_2)} z_{(j_2 i_3)} y_{d,(i_3 j_3)} \nonumber
\\ && \cdots y_{d,(k+m k+m-1)} z_{(k+m-1 k+m)} y_{d,(k+m k+m),}
\quad (d=1,\dots,\hat{s}); \nonumber \\ &&\mbox{ and
the monomial }\quad W(Y,Z)=z_{(11)}\cdot \bigl(
\prod_{d=1}^{\hat{s}} (w_{d} z_{(k+m 1)}) \bigr) \bigr).
\end{eqnarray}
Then the polynomial $f(Y,Z)=\left( \prod_{d=1}^{\hat{s}} (\mathcal{A}_{Y_{d, \bar{0}}} \
\mathcal{A}_{Y_{d,\bar{1}}} )\right) W(Y,Z)$ \
is $(t_{\bar{0}},t_{\bar{1}})$-al\-ter\-na\-ting in any set
$Y_d=Y_{d, \bar{0}} \cup Y_{d, \bar{1}}$ ($d=1,\dots,\hat{s}$).

Notice that for any $d$ and $\theta$ the set $Y_{d, \theta}$ contains at most 1
variable $y_{d,(ij)}$ for the same pair $(i,j).$
Consider the evaluation $y_{d,(ij)}=E_{ij},$ \
$z_{(ij)}=E_{ij}$ ($\ast$) of the polynomial $f.$
Since the variables $z_{(ij)}$ fix the positions for indices
of the variables $y$ then ($\ast$)
gives a nonzero result
\begin{eqnarray} \label{res1}
f_{\bigr|_{(\ast)}}= W_{\bigr|_{(\ast)}} = E_{11} \ne 0.
\end{eqnarray}

Similarly, we construct the polynomial $f(Y,Z),$ and the corresponding nonzero evaluation
for the case $A=M_{k}(F[c]).$ In this case
$A$  has a homogeneous in the grading basis of the type $\{ E_{ij},  E_{ij} c
| i, j = 1,\dots,k \},$ where $E_{ij}$ are the matrix units, and $c$ is the central element
of $A$ satisfying $c^2=1.$ Here $\deg_{\mathbb{Z}_2} E_{ij} =\bar{0},$ and
$\deg_{\mathbb{Z}_2} E_{ij} c =\bar{1}$ for all $i,j=1,\dots,k.$
The set $Y_d=\{y_{d,\bar{0},(ij)}, y_{d,\bar{1},(ij)} \in X^{\mathbb{Z}_2}|
i,j=1,\dots,k \}$ contains $2 k^2$ graded variables. A graded variable $y_{d,\bar{0},(ij)}$
corresponds to the basic element $E_{ij},$ and $y_{d,\bar{1},(ij)}$ corresponds to $E_{ij} \cdot c,$ \
$\deg_{\mathbb{Z}_2} y_{d,\theta,(ij)} = \theta$ ($d=1,\dots,\hat{s}$).
Then $Y_{d, \theta}=\{ y_{d,\theta,(ij)} | i,j=1,\dots,k \},$ and
$|Y_{d, \theta}|=t_{\theta}=\dim A_{\theta}=k^2,$ \ $\theta \in \mathbb{Z}/2 \mathbb{Z}.$
All connecting variables $z_{(ij)}$ have even degree.
In this case the variable $z_{(j_1 i_2)}$ connects variables
$y_{d,\theta,(i_1 j_1)}$ and $y_{d,\xi,(i_2 j_2)}$ for any $\theta, \xi \in \mathbb{Z}/2 \mathbb{Z}.$

The graded monomial $w_{d, \theta}$ is the product of all variables $y_{d,\theta,(i_s j_s)}$
connected by variables $z_{(ij)}$ for any fixed $\theta \in \mathbb{Z}/2 \mathbb{Z}$ and $d$
\begin{eqnarray} \label{word12}
&&w_{d, \theta}= y_{d,\theta,(11)} z_{(11)} y_{d,\theta,(12)} z_{(21)}
y_{d,\theta, (13)} \cdots  y_{d,\theta,(i_1 j_1)} z_{(j_1 i_2)}
y_{d,\theta,(i_2 j_2)} z_{(j_2 i_3)} y_{d,\theta,(i_3 j_3)} \nonumber
\\ && \cdots y_{d,\theta,(k k-1)} z_{(k-1 k)} y_{d,\theta,(k k),}
\quad (d=1,\dots,\hat{s}, \quad \theta = \bar{0}, \bar{1}). \nonumber \\ &&\mbox{
Then the monomial } \quad W(Y,Z)=z_{(11)}\cdot \bigl(
\prod_{d=1}^{\hat{s}} (w_{d \bar{0}} z_{(k 1)} w_{d \bar{1}} z_{(k 1)}) \bigr).
\end{eqnarray}
Take the polynomial $f(Y,Z)=\left( \prod_{d=1}^{\hat{s}}
\mathcal{A}_{Y_{d, \bar{0}}} \cdot \mathcal{A}_{Y_{d, \bar{1}}}\right) W(Y,Z).$ \
$f(Y,Z)$ is $(t_{\bar{0}},t_{\bar{1}})$-al\-ter\-na\-ting in any set
$Y_d=\bigcup_{\theta \in \mathbb{Z}/2 \mathbb{Z}} Y_{d, \theta}$ ($d=1,\dots,\hat{s}$).

For any $d$ and $\theta$ the set $Y_{d, \theta}$ contains at most 1
variable $y_{d,\theta,(ij)}$ for the same pair $(i,j)$ then the
evaluation $y_{d,\bar{0},(ij)}=E_{ij},$ \ $y_{d,\bar{1},(ij)}=E_{ij} c,$
$z_{(ij)}=E_{ij}$ ($\ast$) of $f$
gives a nonzero result
\begin{eqnarray} \label{res12}
f_{\bigr|_{(\ast)}}= W_{\bigr|_{(\ast)}} = E_{11} \cdot c^{k^2} \ne
0.
\end{eqnarray}
Notice that $c^{k^2}=c$ if $k$ is odd, and $c^{k^2}=1$ otherwise.

In both of the cases $f \notin \mathrm{Id}^{\mathbb{Z}_2}(A),$ and this is the desired
polynomial of the type $(t_{\bar{0}},t_{\bar{1}};0;\hat{s}).$ Therefore
$\beta(A)=\mathrm{dims}_{\mathbb{Z}_2} A.$ By Lemma \ref{bound} the
conditions $\dim J(A)=0,$ and $\mathrm{ind}_{\mathbb{Z}_2}(A) =
\mathrm{par}_{\mathbb{Z}_2}(A)$ imply that $A$ is $\mathbb{Z}_2$PI-reduced. \hfill $\Box$

$\mathbb{Z}_2$PI-reduced algebras exist and possess the next properties.

\begin{lemma} \label{reduc}
Given a $\mathbb{Z}_2$PI-reduced algebra with the Wedderburn-Malcev
decomposition (\ref{matrix}) $A=(C_1 \times \cdots \times C_p) \oplus J$
we have $C_{\sigma(1)} J C_{\sigma(2)} J \cdots J C_{\sigma(p)} \ne 0$ for some $\sigma \in
\mathrm{Sym}_p.$
\end{lemma}
\noindent {\bf Proof.}
Suppose that $C_{\sigma(1)} J C_{\sigma(2)} J \cdots J C_{\sigma(p)} = 0$ for any $\sigma \in
\mathrm{Sym}_p.$ Consider the $\mathbb{Z}/2 \mathbb{Z}$-graded subalgebras with elementary decomposition $A_i=(\prod \limits_{\mathop{1 \le j \le p}\limits_{\scriptstyle j \ne i}} C_j )
\oplus J(A)$ of the superalgebra $A.$ Then we have
$\mathrm{Id}^{\mathbb{Z}_2}(A)=\bigcap \limits_{i=1}\limits^{p}
\mathrm{Id}^{\mathbb{Z}_2}(A_i),$ and $\mathrm{dims}_{\mathbb{Z}_2} A_i < \mathrm{dims}_{\mathbb{Z}_2} A$
for any $i=1,\dots,p.$ This contradicts to the definition of $\mathbb{Z}_2$PI-reducible algebra.
\hfill $\Box$

Particularly we have $\mathrm{nd}(A) \ge
p$ for a $\mathbb{Z}_2$PI-reduced algebra $A$.

\begin{lemma} \label{decomp}
Any finite dimensional superalgebra is $\mathbb{Z}_2$PI-equivalent to a
finite direct product of $\mathbb{Z}_2$PI-reduced algebras.
\end{lemma}
\noindent {\bf Proof.} By Lemma \ref{eldec-fd} any finite dimensional superalgebra is $\mathbb{Z}_2$PI-equivalent to a finite dimensional superalgebra $A$ with elementary decomposition.
If $A$ is not $\mathbb{Z}_2$PI-reduced then it is $\mathbb{Z}_2$PI-equivalent to a finite direct
product of finite dimensional elementary decomposed superalgebras with the
complex parameters less than $\mathrm{cpar}_{\mathbb{Z}_2}(A).$ We apply this process
inductively to all multipliers. The set $\mathbb{N}_0^{4}$ with the lexicographical order
satisfies the descending chain condition. Therefore this process of
decomposition will stop after a finite number of steps. \hfill $\Box$

Lemma \ref{decomp} along with Lemmas \ref{subset}, \ref{dirsumKP}
implies

\begin{lemma} \label{Smu}
Any finite dimensional superalgebra $A$ is $\mathbb{Z}_2$PI-equivalent to
a direct product $\mathcal{O}(A) \times \mathcal{Y}(A)$ of finite dimensional superalgebras
with elementary decomposition satisfying
$\mathrm{ind}_{\mathbb{Z}_2}(A)=\mathrm{ind}_{\mathbb{Z}_2}(\mathcal{O}(A))
> \mathrm{ind}_{\mathbb{Z}_2}(\mathcal{Y}(A)).$ Moreover
$\mathcal{O}(A)=A_1 \times \dots \times
A_\rho,$ where $A_i$ are $\mathbb{Z}_2$PI-reduced superalgebras,
and $\mathrm{ind}_{\mathbb{Z}_2}(A_i)=\mathrm{ind}_{\mathbb{Z}_2}(A)$ for all
$i=1,\dots,\rho.$ There exists $\widehat{\mu} \in
\mathbb{N}_0$ such that $S_\mu(A)=S_\mu(\mathcal{O}(A))=\bigcup_{i=1}^\rho S_\mu(A_i)
\subseteq \mathrm{Id}^{\mathbb{Z}_2}(\mathcal{Y}(A))$ holds for any $\mu \ge \widehat{\mu}$.
\end{lemma}
\noindent {\bf Proof.}
Suppose that $A$ is $\mathbb{Z}_2$PI-equivalent to a direct product $A_1 \times \dots
\times A_{\hat{\rho}}$ of $\mathbb{Z}_2$PI-reduced superalgebras.
By Lemma \ref{dirsum} $\mathrm{ind}_{\mathbb{Z}_2}(A)=\max_{1 \le i \le \hat{\rho}} \mathrm{ind}_{\mathbb{Z}_2}(A_i).$ Assume that the Kemer index has the maximal value for $A_i$ with $i=1,\dots,\rho.$ Then the superalgebras $\mathcal{O}(A)=A_1 \times \dots \times A_\rho,$ \
$\mathcal{Y}(A)=A_{\rho+1} \times \dots \times A_{\hat{\rho}}$ satisfy the assertion of the lemma.
\hfill $\Box$

\begin{definition} \label{sen}
$\mathcal{O}(A)$ is called the senior part of $A,$
$\mathcal{Y}(A)$ is called the junior part of $A.$ The
algebras $A_i$ ($i=1,\dots,\rho$) are called the senior components
of $A.$ $\widehat{\mu}(A)$ is the minimal $\widehat{\mu} \in
\mathbb{N}_0$ satisfying the assertion of Lemma \ref{Smu}.
\end{definition}

The next lemma shows the relation between the Kemer index and the parameters
of a $\mathbb{Z}_2$PI-reduced superalgebra.

\begin{lemma} \label{K1}
Given a $\mathbb{Z}_2$PI-reduced superalgebra $A$ we have $\beta(A)=\mathrm{dims}_{\mathbb{Z}_2}
A.$
\end{lemma}
\noindent {\bf Proof.} For a nilpotent superalgebra the assertion is
trivial. Consider any non-nilpotent $\mathbb{Z}_2$PI-reduced algebra $A$ with
$\mathrm{dims}_{\mathbb{Z}_2} A=(t_{\bar{0}},t_{\bar{1}}).$ By Lemma \ref{bound}
it is enough to find a polynomial of the
type $(t_{\bar{0}},t_{\bar{1}};0;\hat{s})$ that is not a graded identity of
$A$ for any $\hat{s} \in \mathbb{N}.$
Assume that $A$ has the elementary decomposition (\ref{matrix}).
For any $l=1,\dots,p$ consider the graded monomial
$W_l(Y_{(l)},Z_{(l)})$ of the type (\ref{word1}) or (\ref{word12}) constructing for the simple component $C_l$ (see Lemma \ref{ind-simple}). $W_l(Y_{(l)},Z_{(l)})$
depends on the disjoint sets of graded variables
$Y_{(l)}=\bigcup_{d=1}^{\hat{s}}  (Y_{(l),(d, \bar{0})} \cup Y_{(l),(d, \bar{1})}),$ and
$Z_{(l)}=\{z_{(l),(ij)}, | i,j=1,\dots,s_l \}.$
Where $Y_{(l)}=\{y_{(l),(d,(ij))} \in X^{\mathbb{Z}_2} \ | \ \ 1 \le i,j \le s_l; \  1 \le d \le
\hat{s} \}$ for $C_l=M_{k_l,m_l}(F),$ \ and $Y_{(l)}=\{y_{(l),(d,\theta,(ij))}
\in X^{\mathbb{Z}_2} \ | \theta = \bar{0}, \bar{1}; \ \ 1 \le i,j \le s_l; \  1 \le d \le
\hat{s} \}$ for  $C_l=M_{k_l}(F[c]).$
We have that $Y_{(l), (d, \theta)} \subseteq X_\theta,$ \ $Y_{d,
\theta}= \bigcup_{l=1}^{p} Y_{(l), (d, \theta)},$ and $\deg_{\mathbb{Z}_2}
z_{(l), (ij)}=\deg_{\mathbb{Z}_2} E_{l i j}$ in $C_l.$ It is clear that $|Y_{d, \theta}|=t_{\theta}$ for any
$d=1,\dots,\hat{s},$ $\theta \in \mathbb{Z}/2 \mathbb{Z}.$

Then the appropriate elementary evaluation of the word $W_l(Y_{(l)},Z_{l})$ in $C_l$ (see Lemma
\ref{ind-simple}) is equal to the non-zero element $E_{l 11} \bar{c}_l,$  where $\bar{c}_l=1$ or $\bar{c}_l=c.$
By Lemma
\ref{reduc} we can assume that $A$ contains an element
$\varepsilon_{1} r_1 \varepsilon_2 \cdots \varepsilon_{p-1}
r_{p-1} \varepsilon_p \ne 0,$ where $r_l \in J$
are some $\mathbb{Z}/2 \mathbb{Z}$-homogeneous radical
elements, and $\varepsilon_l=\sum \limits_{i_l=1}^{s_l} E_{l i_l i_l}$ is the unit of $C_l$
($l=1,\dots,p$).
Let us take $p$ nontrivial graded polynomials
\begin{equation} \label{polynK1}
f_l(\widetilde{X}_{(l)},Y_{(l)},Z_{(l)})=\sum_{i_l=1}^{s_l} \tilde{x}_{l,(i_l)}
W_l(Y_{(l)},Z_{(l)}) \ \tilde{\tilde{x}}_{l,(i_l,\bar{c}_l)}
\end{equation}
depending on the additional set of graded variables
$\widetilde{X}_{(l)}=\{
\tilde{x}_{l,(i_l)},  \tilde{\tilde{x}}_{l,(i_l,
\bar{c}_l)} \in X^{\mathbb{Z}_2} \  | \ 1 \le
i_l \le s_l; \deg_{\mathbb{Z}_2} \tilde{x}_{l,(i_l)} =
\deg_{\mathbb{Z}_2} E_{l i_l 1}, \
\deg_{\mathbb{Z}_2} \tilde{\tilde{x}}_{l,(i_l,\bar{c}_l)} =
\deg_{\mathbb{Z}_2} E_{l 1 i_l} \bar{c}_l,\}.$
Notice that $f_l$ is not a multihomogeneous polynomial,
although it is linear in $Y_{(l)},$ and $Z_{(l)}.$
But $f_l$ is homogeneous in the grading of even degree,
since any monomial $\tilde{x}_{l,(i_l)}
W_l(Y_{(l)},Z_{(l)}) \ \tilde{\tilde{x}}_{l,(i_l,\bar{c}_l)}$
has even degree (it has the same $\mathbb{Z}/2 \mathbb{Z}$-degree as
$E_{l i_l i_l}$ in $C_l$).

Then the polynomial
\begin{eqnarray} \label{polyn2}
f(\widetilde{X},Y,Z)=\left( \prod_{d=1}^{\hat{s}}
\mathcal{A}_{Y_{d, \bar{0}}} \  \mathcal{A}_{Y_{d, \bar{1}}} \right) (f_1 x_1 f_2 x_2 \cdots
x_{p-1} f_p)
\end{eqnarray}
is linear in $Y \bigcup Z=
\left( \bigcup_{l=1}^{p} Y_{(l)} \right) \bigcup \left(
\bigcup_{l=1}^{p} Z_{(l)} \right),$ and alternating in any set $Y_{d, \theta}
\subseteq X_{\theta},$  $\theta \in \mathbb{Z}/2 \mathbb{Z},$  $d=1,\dots,s.$
It follows from Lemma \ref{ind-simple} that the evaluation
\begin{eqnarray*}
&&y_{(l), (d,(i_l j_l))}=\varepsilon_l E_{l
i_l j_l} \varepsilon_l, \qquad \qquad \
y_{(l), (d,\bar{0},(i_l j_l))}=\varepsilon_l E_{l
i_l j_l} \varepsilon_l, \\
&&y_{(l), (d,\bar{1},(i_l j_l))}=\varepsilon_l (E_{l
i_l j_l} c) \varepsilon_l, \qquad \
z_{(l),(i_l
j_l)}=\varepsilon_l E_{l i_l j_l} \varepsilon_l, \\
&&\tilde{x}_{l,(i_l)}=\varepsilon_l E_{l i_l 1} \varepsilon_l, \qquad \qquad \qquad \ \
\tilde{\tilde{x}}_{l,(i_l,\bar{c}_l)}=
\varepsilon_l (E_{l 1 i_l} \bar{c}_l) \varepsilon_l, \\
&&x_{q}=\varepsilon_{q} r_{q} \varepsilon_{q+1}, \\
&&1 \le i_l, j_l \le
s_l, \  \
l=1,\dots,p; \ \
q=1,\dots, p-1; \ \
d=1,\dots,\hat{s},
\end{eqnarray*}
of the polynomial $f(\widetilde{X},Y,Z)$
is equal to $\varepsilon_{1} r_1 \varepsilon_2 \cdots \varepsilon_{p-1}
r_{p-1} \varepsilon_p \ne 0.$

Therefore we obtain $f \notin \mathrm{Id}^{\mathbb{Z}_2}(A).$ Hence at least one
multihomogeneous component $\tilde{f}$ of $f$ also is not
a graded identity of $A.$ $\tilde{f}$ is alternating in any set $Y_{d,
\theta}$ ($\theta \in \mathbb{Z}/2 \mathbb{Z},$ \ $d=1,\dots,\hat{s}$). Thus
$\tilde{f}$ is the required polynomial. \hfill $\Box$

\section{Exact polynomials.}

\begin{definition}
Given a finite dimensional superalgebra $A$ with elementary decomposition
an elementary evaluation $(a_1,\dots,a_n) \in A^n$ (namely, $a_i \in D \bigcup U
\subseteq A$ (\ref{basisD}), (\ref{basisU})) is called incomplete if there exists
$j=1,\dots,p$ such that $$\{a_1,\dots,a_n\} \bigcap (C_j \bigcup \{
\varepsilon_j u \varepsilon_l, \varepsilon_l u \varepsilon_j | u
\in U, \  l=1, \dots, p+1\}) = \emptyset.$$ Otherwise $(a_1,\dots,a_n)$ is called complete.
\end{definition}

\begin{definition}
An elementary evaluation $(a_1,\dots,a_n) \in A^n$ is
called thin if it contains less than $\mathrm{nd}(A)-1$
radical elements (not necessarily distinct).
\end{definition}

\begin{definition}
We say that a multilinear graded polynomial $f(x_1,\dots,x_n) \in F\langle
X^{\mathbb{Z}_2} \rangle$ is exact for a finite dimensional
superalgebra $A$ with elementary decomposition if $f(a_1,\dots,a_n)=0$ holds in $A$ for
any thin or incomplete evaluation $(a_1,\dots,a_n) \in A^n.$
\end{definition}

\begin{lemma} \label{Exact2}
If $A$ is a $\mathbb{Z}_2$PI-reduced superalgebra then
any multilinear polynomial of the type $(\mathrm{dims}_{\mathbb{Z}_2}
A;\mathrm{nd}(A)-1;0)$ is exact for $A.$
\end{lemma}
\noindent {\bf Proof.} It is clear that a $\mathbb{Z}_2$PI-reduced algebra
either is not semisimple or is $\mathbb{Z}/2 \mathbb{Z}$-simple.
For a $\mathbb{Z}/2 \mathbb{Z}$-graded simple
finite dimensional algebra any multilinear graded polynomial can be
assumed exact. A multilinear polynomial of the type
$(0,0;\mathrm{nd}(A)-1;0)$ has degree greater or equal to
$(\mathrm{nd}(A)-1).$ Hence it is assumed to be exact for a
nilpotent algebra $A.$

Suppose that $A$ is not nilpotent, and $f$ is a multilinear
polynomial of the type $(\mathrm{dims}_{\mathbb{Z}_2} A;\mathrm{nd}(A)-1;0).$
In a thin evaluation at least one of
$\tilde{s}=\mathrm{nd}(A)-1$ collections of $\tau_{j}$-alternating
variables of $f$ will be completely replaced by semisimple
elements. Since $\tau_{j}>\mathrm{dims}_{\mathbb{Z}_2} A$ for any
$j=1,\dots,\tilde{s}$ then the result of such evaluation will be zero.

Any simple $\mathbb{Z}/2 \mathbb{Z}$-graded component $C_l$ of $A$
(Lemma \ref{Pierce}) has the unit $\varepsilon_l \in C_{l \bar{0}}.$
Hence $\dim_F C_{l \bar{0}} > 0$ for any $l=1,\dots,p$. Therefore an incomplete evaluation
can not contain all semisimple elements of the even degree from the base
$D$ (\ref{basisD}). Taking into account the conditions $\tau_{j}>\mathrm{dims}_{\mathbb{Z}_2} A$
for any $j=1,\dots,\tilde{s}$ we obtain that at least two variables of every collection
of $\tau_j$-alternating variables of $f$ must be substituted by
radical elements, otherwise the result of the substitution will be
zero. Thus in any case the result of an incomplete evaluation of the
polynomial $f$ is zero. \hfill $\Box$

\begin{lemma} \label{Exact1}
Any nonzero $\mathbb{Z}_2$PI-reduced superalgebra $A$ has an exact polynomial, that
is not a graded identity of $A.$
\end{lemma}
\noindent {\bf Proof.} For a nilpotent superalgebra $A$ the assertion
follows from Lemma \ref{Exact2}.
Suppose that $A$ is a non-nilpotent superalgebra
with the decomposition (\ref{matrix}). Consider its subalgebras
$A_i=(\prod \limits_{\mathop{1 \le j \le p}\limits_{\scriptstyle j
\ne i}} C_j ) \oplus J(A),$ \ $i=1,\dots, p.$ Take $q=\dim_F
J(A),$ \ $s=\mathrm{nd}(A)-1.$ Then by Lemma \ref{Aqs} \
$\widetilde{A}=A_1 \times \dots \times A_p \times
\mathcal{R}_{q,s}(A)$ is a finite dimensional superalgebra
with elementary decomposition
satisfying $\mathrm{Id}^{\mathbb{Z}_2}(A) \subseteq
\mathrm{Id}^{\mathbb{Z}_2}(\widetilde{A}).$
Let $\{r_1,\dots,r_q\} \subseteq U$ be a
homogeneous basis of $J(A)$ (\ref{basisU}), \ $\deg_{\mathbb{Z}_2} r_i=\theta_i$
($i=1,\dots,q$). Consider the map $\varphi$ defined by the next equalities
$\varphi(x_{i \theta_i})=r_i$ for $i=1,\dots,q$, and $\varphi(b)=b$
for any $b \in B.$ $\varphi$ can be extended to a surjective graded
homomorphism $\varphi: B(X_q^{\mathbb{Z}_2}) \rightarrow A.$
It follows that any multilinear polynomial $f \in \mathrm{Id}^{\mathbb{Z}_2}(\mathcal{R}_{q,s}(A))$
is turned into zero under any thin substitution. It is also clear that any
incomplete substitution in a multilinear graded polynomial
$f \in \mathrm{Id}^{\mathbb{Z}_2}(\times_{i=1}^p A_i)$ yields zero.
Therefore, any multilinear polynomial $f \in \mathrm{Id}^{\mathbb{Z}_2}(\widetilde{A})$ is exact
for $A.$ Remark that $\mathrm{cpar}_{\mathbb{Z}_2}(A_i) < \mathrm{cpar}_{\mathbb{Z}_2}(A)$ ($1 \le i \le
p$), and $\mathrm{cpar}_{\mathbb{Z}_2}(\mathcal{R}_{q,s}(A)) < \mathrm{cpar}_{\mathbb{Z}_2}(A).$
Since $A$ is $\mathbb{Z}_2$PI-reduced then $\mathrm{Id}^{\mathbb{Z}_2}(A) \subsetneqq \mathrm{Id}^{\mathbb{Z}_2}(\widetilde{A}).$ Any multilinear graded polynomial $f$ such that
$f \in \mathrm{Id}^{\mathbb{Z}_2}(\widetilde{A}),$ and $f \notin \mathrm{Id}^{\mathbb{Z}_2}(A)$
satisfies the assertion of the lemma.  \hfill
$\Box$

\begin{lemma} \label{Gammasub}
Let $A$ be a finite dimensional superalgebra with an elementary decomposition, $h$ an
exact polynomial for $A,$ and $\bar{a} \in A^n$ is any complete
evaluation of $h$ containing $\tilde{s}=\mathrm{nd}(A)-1$ radical elements.
Then for any $\mu
\in \mathbb{N}_0$ there exist a graded polynomial
$h_\mu \in \mathbb{Z}_2 T[h],$ and
an elementary evaluation $\bar{u}$ of $h_\mu$ by elements of
$A$ such that:
\begin{enumerate}
 \item $h_\mu(\widetilde{Y}_1,\dots,\widetilde{Y}_{\tilde{s}+\mu},
\widetilde{X},\widetilde{Z})$ is
$\tau_j$-alternating in any set $\widetilde{Y}_j$ with $\tau_j >
\beta=\mathrm{dims}_{\mathbb{Z}_2} A$ for all $j=1,\dots,\tilde{s}$, and
is $\beta$-alternating in any $\widetilde{Y}_j$ for
$j=\tilde{s}+1,\dots,\tilde{s}+\mu$ (all the sets
$\widetilde{Y}_j,$ $\widetilde{X},$ $\widetilde{Z}$ are disjoint),
 \item $h_\mu(\bar{u})=h(\bar{a}),$
 \item all the
variables of $\widetilde{X} \bigcup \widetilde{Z}$ are
replaced by semisimple elements.
\end{enumerate}
\end{lemma}
\noindent {\bf Proof.}
If $h(\bar{a})=0$ then the assertion of lemma is trivial. It is sufficient to take any consequence
$h_\mu(\widetilde{Y}_1,\dots,\widetilde{Y}_{\tilde{s}+\mu}) \in \mathbb{Z}_2 T[h]$
that is alternating in all $\widetilde{Y}_j$ as required (we assume here that $\widetilde{X} \bigcup \widetilde{Z}=\emptyset$), and replace the variables of the alternating sets $\widetilde{Y}_j$  by equal elements. Particularly, from conditions $\mathrm{nd}(A)=1,$ \ $p \ge 2$ it follows that $h(\bar{a}) = 0.$

Assume that  $h(\bar{a}) \ne 0.$
Consider the case $\mathrm{nd}(A)>1,$ \ $p
\ge 2$ in the decomposition (\ref{matrix}) of the superalgebra $A.$ We can assume
for simplicity that the evaluation $\bar{a}$ has the form
\quad $\bar{a}=(\mathbf{\varepsilon_1} r_{\mathbf{1}}
\varepsilon_{l''_1},\varepsilon_{l'_2} r_{\mathbf{2}}
\mathbf{\varepsilon_2},\dots,\varepsilon_{\mathbf{p-1}} r_{\mathbf{q}} \mathbf{\varepsilon_{p}},$
$\varepsilon_{l'_{q+1}} r_{q+1} \varepsilon_{l''_{q+1}},\dots,$
$\varepsilon_{l'_{\tilde{s}}} r_{\tilde{s}}
\varepsilon_{l''_{\tilde{s}}},b_1,\dots,b_{n-\tilde{s}}),$ \ where
$1 \le l'_s,l''_s \le p+1,$ \ $\{1,\dots,p\} \subseteq \{l'_s, l''_s
\ | 1 \le s \le q \}$ for some $q \le \tilde{s},$ \ $l'_s \ne l''_s$ for any $1 \le s \le q.$
Namely, in a complete evaluation all minimal orthogonal idempotents $\varepsilon_1,\dots,\varepsilon_p$
appear in mixed radical elements (in elements of the type $\varepsilon_{l'} r_l \varepsilon_{l''} \in U$ \ (\ref{basisU}) with $l' \ne l''$). We assume that they appear in the first $q$ mixed radical elements.
The last elements $n-\tilde{s}$ of the evaluation $\bar{a}$ \ $b_1,\dots,$ $b_{n-\tilde{s}} \in D$ (\ref{basisD}) are supposed to be semisimple, and the first $\tilde{s}$ elements
$\varepsilon_{l'_{s}} r_{s} \varepsilon_{l''_{s}} \in U$ \ (\ref{basisU}) are radical.

Let us take for any $l=1,\dots,p$ the polynomial
$f_l(\widetilde{X}_{(l)},Y_{(l)},Z_{(l)})$ defined by
(\ref{polynK1}) in Lemma \ref{K1}  assuming
$\hat{s}=\tilde{s}+\mu=\mathrm{nd}(A)-1+\mu.$ Here
$Y_{(l)}=\bigcup_{d=1}^{\tilde{s}+\mu} (Y_{(l),(d, \bar{0})} \cup
Y_{(l),(d, \bar{1})}),$  and $Y_{d, \theta}= \bigcup_{l=1}^{p} (Y_{(l),(d, \bar{0})} \cup
Y_{(l),(d, \bar{1})})$ ($d=1,\dots,\tilde{s}+\mu,$ \ $\theta \in \mathbb{Z}/2 \mathbb{Z}$).
consider also a new set of graded variables $\tilde{y}_s \in X^{\mathbb{Z}_2}$
such that $\deg_{\mathbb{Z}_2} \
\tilde{y}_s =\deg_{\mathbb{Z}_2} \ r_s,$ \ $1 \le s \le \tilde{s}.$ Let us denote $\widetilde{Y}_{d,
\theta}=Y_{d, \theta} \bigcup \{ \tilde{y}_d \}$ if $\deg_{\mathbb{Z}_2} \
\tilde{y}_d =\theta,$ and $\widetilde{Y}_{d,
\theta}=Y_{d, \theta}$ otherwise $d=1,\dots,\mu+\tilde{s}$).

Then we obtain $\bar{\tau}_d=(\tau_{d \bar{0}},\tau_{d \bar{1}})=
(|\widetilde{Y}_{d, \bar{0}}|,|\widetilde{Y}_{d,\bar{1}}|) >
(|Y_{d, \bar{0}}|,|Y_{d, \bar{1}}|)$ if $d=1,\dots,\tilde{s}.$  And
$\bar{\tau}_d=(|\widetilde{Y}_{d,\bar{0}}|,|\widetilde{Y}_{d, \bar{1}}|) =
(|Y_{d, \bar{0}}|,|Y_{d, \bar{1}}|)$ for any
$d=\tilde{s}+1,\dots,\tilde{s}+\mu.$ Here
$(|Y_{d,\bar{0}}|,|Y_{d, \bar{1}}|)=\mathrm{dims}_{\mathbb{Z}_2} A$ for any $d.$

Denote by $\widetilde{Y}_d=\widetilde{Y}_{d, \bar{0}} \cup \widetilde{Y}_{d, \bar{1}}$ \ ($d=1,\dots,\tilde{s}+\mu$); \
$\widetilde{X}=\bigl(\bigcup_{l=1}^{p} \widetilde{X}_l \bigr)
\bigcup \{x_1,\dots,x_{n-\tilde{s}}\};$ \
$\widetilde{Z}=\bigcup_{l=1}^{p} Z_l.$ Consider the
polynomials
\begin{eqnarray*}
&&h'(\widetilde{Y}_1,\dots,\widetilde{Y}_{\tilde{s}+\mu},\widetilde{X},\widetilde{Z})=
h \Bigl(f_1 \bigl(\widetilde{X}_{(1)},Y_{(1)},Z_{(1)} \bigr)
\tilde{y}_1, \ \tilde{y}_2 f_2
\bigl(\widetilde{X}_{(2)},Y_{(2)},Z_{(2)} \bigr), \ \dots, \\
&&f_{p-1}\bigl(\widetilde{X}_{(p-1)},Y_{(p-1)},Z_{(p-1)} \bigr)
\tilde{y}_q f_p\bigl(\widetilde{X}_{(p)},Y_{(p)},Z_{(p)} \bigr), \
\tilde{y}_{q+1},\dots,\tilde{y}_{\tilde{s}}, \
x_1,\dots,x_{n-\tilde{s}} \Bigr); \\
&&h_\mu(\widetilde{Y}_1,\dots,\widetilde{Y}_{\tilde{s}+\mu},
\widetilde{X},\widetilde{Z})= \Bigl( \prod_{d=1}^{\tilde{s}+\mu} (\mathcal{A}_{\widetilde{Y}_{d, \bar{0}}} \
\mathcal{A}_{\widetilde{Y}_{d, \bar{1}}} ) \Bigr)
h'(\widetilde{Y}_1,\dots,\widetilde{Y}_{\tilde{s}+\mu},\widetilde{X},\widetilde{Z}).
\end{eqnarray*}
The polynomial $h_\mu$ is linear in variables
$\widetilde{Y}=\bigcup_{d=1}^{\tilde{s}+\mu} \widetilde{Y}_d,$ and
$\bar{\tau}_d$-alternating in any $\widetilde{Y}_d$ \
($d=1,\dots,\tilde{s}+\mu$). Although $h_\mu$ is not multilinear and
is not multihomogeneous in general.

Consider the following evaluation of the polynomial $h_\mu$ in the superalgebra $A$
\begin{eqnarray} \label{substK22}
&&y_{(l), (d,(i_l j_l))}=\varepsilon_l E_{l i_l j_l} \varepsilon_l, \qquad \qquad \
y_{(l), (d,\bar{0},(i_l j_l))}=\varepsilon_l E_{l i_l j_l} \varepsilon_l, \nonumber\\
&&y_{(l), (d,\bar{1},(i_l j_l))}=\varepsilon_l (E_{l i_l j_l} c) \varepsilon_l, \qquad \
z_{(l),(i_l j_l)}=\varepsilon_l E_{l i_l j_l} \varepsilon_l, \nonumber \\
&&\tilde{x}_{l,(i_l)}=\varepsilon_l E_{l i_l 1} \varepsilon_l, \qquad \qquad \qquad \ \
\tilde{\tilde{x}}_{l,(j_l,\bar{c}_l)}=
\varepsilon_l (E_{l 1 j_l} \bar{c}_l) \varepsilon_l, \nonumber \\
&&\tilde{y}_s=a_s=\varepsilon_{l'_s} r_s \varepsilon_{l''_s}; \qquad \qquad \qquad x_{n'}=
a_{n'+\tilde{s}}=b_{n'}, \\
&&l=1,\dots,p; \qquad \qquad \qquad 1 \le i_l, j_l \le s_l,  \nonumber\\
&&d=1,\dots,\tilde{s}+\mu; \qquad 1 \le s \le \tilde{s}; \qquad 1
\le n' \le n-\tilde{s}. \nonumber
\end{eqnarray}
The evaluations of the variables $y,$ $z,$ $\tilde{x},$ and
$\tilde{\tilde{x}}$ are defined as in Lemma \ref{K1}. It is clear that (\ref{substK22}) is
an elementary evaluation and satisfies the third claim of the lemma.

Due to the evaluation of the variables $z,$ $\tilde{x},$ and
$\tilde{\tilde{x}}$, the polynomial $f_l$ can contain only elements of
the simple component $C_l$ or elements $\varepsilon_l r_s \varepsilon_l$, otherwise, we get zero.
The second case will give us a thin evaluation of the polynomial $h,$ thus, such summands are also zero. Therefore, the evaluation (\ref{substK22})
of the polynomial $h_\mu$ gives the same result as this evaluation of the polynomial
\[ h(f'_1 \tilde{y}_1, \ \tilde{y}_2 f'_2, \ \dots,
f'_{p-1} \tilde{y}_q f'_p, \  \tilde{y}_{q+1},\dots,\tilde{y}_{\tilde{s}}, \
x_1,\dots,x_{n-\tilde{s}}),\]
where $f'_l=( \prod_{d=1}^{\tilde{s}+\mu} \mathcal{A}_{\widetilde{Y}_{d, \bar{0}}} \
\mathcal{A}_{\widetilde{Y}_{d, \bar{1}}} ) f_l.$
Similarly to arguments of Lemmas \ref{K1}, \ref{ind-simple},
we can see that the result of our evaluation of $f'_l$ is equal to $\varepsilon_l.$
Thus (\ref{substK22}) for $h_\mu$ gives the result
$h(\varepsilon_1 r_1 \varepsilon_{l''_1}, \ \varepsilon_{l'_2}
r_2 \varepsilon_2, \ \dots, \ \varepsilon_{p-1} r_{q}
\varepsilon_{p},\ \varepsilon_{l'_{q+1}} r_{q+1} \varepsilon_{l''_{q+1}},\dots,
\varepsilon_{l'_{\tilde{s}}} r_{\tilde{s}} \varepsilon_{l''_{\tilde{s}}},b_1,\dots,b_{n-\tilde{s}})=
h(a_1,\dots,a_n).$ Therefore $h_\mu$ is the desired polynomial. The evaluation (\ref{substK22})
can be taken as $\bar{u}$.

Consider the case $\mathrm{nd}(A) \ge 1,$ and $p=1.$
If $\varepsilon_1 h(a_1,\dots,a_n)=0$ then among
$\tilde{s}$ radical elements of $\bar{a}$ there is a
homogeneous in the grading radical element of the type
$\varepsilon_2 r \varepsilon_1 \in \varepsilon_2 J(A)
\varepsilon_1,$ where $\varepsilon_2$ is the adjoint idempotent of
$A$ (since the substitution $\bar{a}$ is complete).
We can suppose
without loss of generality that $a_1=\varepsilon_2 r_1
\varepsilon_1,$ and the first $\tilde{s}=(\mathrm{nd}(A)-1)$
elements of $\bar{a}$ are radical $a_2=\varepsilon_{l'_2} r_2
\varepsilon_{l''_2},\dots,
a_{\tilde{s}}=\varepsilon_{l'_{\tilde{s}}} r_{\tilde{s}}
\varepsilon_{l''_{\tilde{s}}} \in J(A),$ $a_{n'+\tilde{s}}=b_{n'}
\in D,$ for $n'=1,\dots,n-\tilde{s}.$ Then using the same arguments
as in the previous case $p \ge 2,$ assuming that $p=1,$ $q=1,$ \
$l'_{s}, l''_{s} \in \{1, 2\},$ we can prove that the evaluation
(\ref{substK22}) of the polynomial
\[ h_\mu=
\Bigl( \prod_{d=1}^{\tilde{s}+\mu} \mathcal{A}_{\widetilde{Y}_{d, \bar{0}}}
\mathcal{A}_{\widetilde{Y}_{d, \bar{1}}} \Bigr) \  h
\Bigl(\tilde{y}_1 \cdot f_1
\bigl(\widetilde{X}_{(1)},Y_{(1)},Z_{(1)} \bigr), \
\tilde{y}_2,\dots,\tilde{y}_{\tilde{s}}, x_1,\dots,x_{n-\tilde{s}} \Bigr)
\]
is equal to $h(a_1,\dots,a_n).$

If $\varepsilon_1  h(a_1,\dots,a_n)\ne 0$ then the similar arguments
as in the case $p \ge 2$ show that the result of the
evaluation (\ref{substK22}) of the polynomial
\begin{equation*}
h_\mu= \Bigl( \prod_{d=1}^{\tilde{s}+\mu} (\mathcal{A}_{\widetilde{Y}_{d, \bar{0}}} \
\mathcal{A}_{\widetilde{Y}_{d, \bar{1}}}) \Bigr) \ f_1
\bigl(\widetilde{X}_{(1)},Y_{(1)},Z_{(1)} \bigr) \cdot h
\Bigl(\tilde{y}_1,\tilde{y}_2,\dots,\tilde{y}_{\tilde{s}},
x_1,\dots,x_{n-\tilde{s}} \Bigr)
\end{equation*}
is equal to $\varepsilon_1 h(a_1,\dots,a_n)= h(a_1,\dots,a_n).$
Here variables $\tilde{y}_s$ can not take places
inside $f_1$ with a nonzero result,
since it will give a thin evaluation of $h.$

In both of the last cases the polynomial $h_\mu$ and the corresponding
evaluation possess all desired properties.

The case $p=0$ is trivial. In this case the superalgebra $A$ is nilpotent, and
$\mathrm{dims}_{\mathbb{Z}_2}  A=(0,0).$ Since a multihomogeneous graded
polynomial $h$ is exact for $A$ and is not a graded identity of $A$ ($h(\bar{a}) \ne 0$)
then the full degree of $h$ is $(\mathrm{nd}(A)-1),$ and $h$ has the
type $(0,0;\mathrm{nd}(A)-1;\mu)=(\mathrm{dims}_G
A;\mathrm{nd}(A)-1;\mu)$ for any $\mu \in \mathbb{N}_0.$ Thus in this case
$h_\mu=h,$ \ $\bar{u}=\bar{a}.$
\hfill $\Box$

Lemma \ref{Gammasub} has the following important corollaries.

\begin{lemma} \label{gamma}
Let $A$ be a $\mathbb{Z}_2$PI-reduced algebra then \
$\mathrm{ind}_{\mathbb{Z}_2}(A)=\mathrm{par}_{\mathbb{Z}_2}(A).$ If $f$ is an exact polynomial
for $A$, and $f \notin \mathrm{Id}^{\mathbb{Z}_2}(A)$ then $\mathbb{Z}_2 T[f] \bigcap
S_{\mu}(A) \ne \emptyset$ for any $\mu \in \mathbb{N}_0.$
\end{lemma}
\noindent {\bf Proof.} By Lemma \ref{Exact1} $A$ has an exact
polynomial $\tilde{f} \notin \mathrm{Id}^{\mathbb{Z}_2}(A).$ $\tilde{f}$ can be
nonzero only for a complete evaluation containing exactly
$(\mathrm{nd}(A)-1)$ radical elements. Lemma \ref{Gammasub}
implies that $\tilde{f}$ has a nontrivial consequence $\widetilde{g}
\notin \mathrm{Id}^{\mathbb{Z}_2}(A)$ of the type $(\mathrm{dims}_{\mathbb{Z}_2}
A;\mathrm{nd}(A)-1;\mu)$ for any $\mu \in \mathbb{N}_0.$ By
Lemma \ref{K1} we have $\beta(A)=\mathrm{dims}_{\mathbb{Z}_2} A.$ Then Definition
\ref{defgamma} implies $\gamma(A) > \mathrm{nd}(A)-1.$ Taking into
account that $\mathrm{ind}_{\mathbb{Z}_2}(A) \le \mathrm{par}_{\mathbb{Z}_2}(A)$ we obtain
$\gamma(A)=\mathrm{nd}(A).$

Moreover by Lemma \ref{Gammasub} any exact for $A$ polynomial $f$ such that
$f \notin \mathrm{Id}^{\mathbb{Z}_2}(A)$  has a nontrivial consequence $g_\mu \in
\mathbb{Z}_2 T[f]$ for any $\mu \in
\mathbb{N}_0.$ Where $g_\mu$ is $\mu$-boundary polynomial for $A.$
\hfill $\Box$

Lemma \ref{gamma} along with Lemma \ref{Exact2} immediately
implies

\begin{lemma} \label{Exact3}
Any multilinear $\mu$-boundary polynomial for a $\mathbb{Z}_2$PI-reduced
algebra $A$ is exact for $A$.
\end{lemma}

\begin{lemma} \label{Exact4}
Given a $\mathbb{Z}_2$PI-reduced algebra $A,$ and an integer $\mu \in
\mathbb{N}_0$ let $S_{A,\mu}$ be any set of graded polynomials of the
type $(\beta(A);\gamma(A)-1;\mu).$ Then if a multilinear polynomial
$f \in \mathbb{Z}_2 T[S_{A,\mu}]+\mathrm{Id}^{\mathbb{Z}_2}(A)$ then $f$ is exact for $A.$
\end{lemma}
\noindent {\bf Proof.} A multilinear graded polynomial $f \in \mathbb{Z}_2 T[S_{A,\mu}]+\mathrm{Id}^{\mathbb{Z}_2}(A)$ has the
form
\[
f(x_1,\dots,x_n)=\sum \limits_{\mathop{1 \le i \le
d_1}\limits_{\scriptstyle 1 \le j \le d_2}} \alpha_{i j} \ \  v_i
\cdot g_j(h_{i 1},\dots,h_{i n_j}) \cdot w_i + g(x_1,\dots,x_n),
\]
where $h_{i 1},\dots,h_{i
n_j},$ $v_i,$ $w_i \in F\langle X^{\mathbb{Z}_2} \rangle$ are graded monomials,
$v_i,$ $w_i$ are possibly empty;
\ $g_j(z_{j 1},\dots,z_{j n_j}) \in F\langle X^{\mathbb{Z}_2} \rangle$ are full
linearizations of some polynomials from the set $S_{A,\mu};$
$g$ is a multilinear graded identity of $A;$ $\alpha_{i j} \in F.$
It is clear that the multilinear graded polynomials $g_j(z_{j
1},\dots,z_{j n_j})$ have the type
$(\beta(A);\gamma(A)-1;\mu)$  ($\forall j=1,\dots,d_2$), and the
graded monomials $v_i \cdot h_{i 1} \cdots h_{i n_j}
\cdot w_i$ are multilinear and depends on the same
variables $x_1,\dots,x_n$ as $f$ for any $i.$
Particularly, the monomials $h_{i
1}(x_{\delta_1},\dots,x_{\delta_{s_{i 1}}}),$ $\dots, h_{i
n_j}(x_{\lambda_1},\dots,x_{\lambda_{s_{i n_j}}})$ are multilinear, here we have
$\{ \delta_1,\dots,\delta_{s_{i 1}},\dots,
\lambda_1,\dots,\lambda_{s_{i n_j}} \} \subseteq \{ 1, \dots, n\}.$

Fix any $i=1,\dots d_1,$ $j=1,\dots,d_2.$
Given an elementary evaluation $(a_1,\dots,a_n)$ of $f$ in $A$  consider
homogeneous elements of $A,$ that are the evaluation of the monomials
$h_{i l}$ \ \  $\widetilde{a}_{i 1}=h_{i 1}(a_{\delta_1},\dots,a_{\delta_{s_{i
1}}}),\dots, \widetilde{a}_{i n_j}=h_{i
n_j}(a_{\lambda_1},\dots,a_{\lambda_{s_{i n_j}}}).$
Then
$\widetilde{a}_{i l}$ is either a semisimple element of $A$ or radical, and
\begin{eqnarray*}
&&\widetilde{a}_{i l}= \sum_{t_l=1}^{\dim \mathcal{R}_l} \alpha_{(i l), t_l}
\ u_{t_l}, \qquad \ \quad \alpha_{(i l), t_l} \in F, \quad
\deg_{\mathbb{Z}_2} u_{t_l} = \deg_{\mathbb{Z}_2} \widetilde{a}_{i l}=\theta_l, \ \forall t_l;\\
&&u_{t_l} \in D \ (\ref{basisD}),  \ \forall t_l, \ \  \mathcal{R}_l=B_{\theta_l}, \qquad
\mbox{ or } \ \
 \qquad u_{t_l} \in U  \ (\ref{basisU}), \  \forall t_l, \ \
\mathcal{R}_l=J(A)_{\theta_l}.
\end{eqnarray*}
Then the evaluation $(a_1,\dots,a_n)$ of $g_j(h_{i
1},\dots,h_{i n_j})$ yields
\begin{eqnarray*}
&&g_j(\widetilde{a}_{i 1},\dots,\widetilde{a}_{i n_j})=
g_j(\sum_{t_1=1}^{\dim \mathcal{R}_1} \alpha_{(i 1), t_1}
u_{t_1},\dots,\sum_{t_{n_j}=1}^{\dim \mathcal{R}_{n_j}} \alpha_{(i n_j),
t_{n_j}} u_{t_{n_j}})= \nonumber \\
&&\sum_{t_1,\dots,t_{n_j}} \alpha_{(i 1), t_1} \cdots \alpha_{(i
n_j), t_{n_j}} \  g_j(u_{t_1},\dots,u_{t_{n_j}}).
\end{eqnarray*}
All evaluations $(u_{t_1},\dots,u_{t_{n_j}})$ of the polynomials
$g_j$ are elementary. The numbers of radical elements in
$(u_{t_1},\dots,u_{t_{n_j}})$ is equal to the number of
radical elements in $(\widetilde{a}_{i 1},\dots,\widetilde{a}_{i
n_j})$ for all $t_1,\dots,t_{n_j},$ and does not exceed the number
of radical elements in the initial evaluation $(a_1,\dots,a_n).$
Thus if $(a_1,\dots,a_n)$ is a thin evaluation then
$(u_{t_1},\dots,u_{t_{n_j}})$ is also thin for all
$t_1,\dots,t_{n_j}$ and $i,$ $j.$

$(\prod \limits_{\mathop{1 \le l \le p}\limits_{\scriptstyle l
\ne k}} C_l ) \oplus (\sum \limits_{\mathop{1 \le l',l'' \le
p+1}\limits_{\scriptstyle l' \ne k, \ l'' \ne k}} \varepsilon_{l'}
J(A) \varepsilon_{l''} )$ is a graded subalgebra of $A$ for any $k.$
Therefore, if $(a_1,\dots,a_n)$ is incomplete
then $(\widetilde{a}_{i 1},\dots,\widetilde{a}_{i n_j}),$ and $(u_{t_1},\dots,u_{t_{n_j}})$
are also incomplete. The last evaluations do not contain elements of
$C_k \bigcup \{ \varepsilon_k r
\varepsilon_l, \varepsilon_l r \varepsilon_k | r \in U, \  1 \le l
\le p+1 \}$ if elements of this set does not appear in $(a_1,\dots,a_n).$

$A$ is $\mathbb{Z}_2$PI-reduced, and by Lemma
\ref{gamma} $\mathrm{ind}_{\mathbb{Z}_2}(A)=\mathrm{par}_{\mathbb{Z}_2}(A).$
Thus $g_j$ is a multilinear graded polynomial of
the type $(\beta(A);\gamma(A)-1;\mu)=(\mathrm{dims}_{\mathbb{Z}_2}
A;\mathrm{nd}(A)-1;\mu),$ and  by Lemma \ref{Exact2} $g_j$ is exact for $A$ (for any $j$).

Thus for any thin or incomplete evaluation $(a_1,\dots,a_n)$ we obtain
\[g_j(\widetilde{a}_{i 1},\dots,\widetilde{a}_{i
n_j})=\sum_{t_1,\dots,t_{n_j}} \alpha_{(i l), t_1} \cdots \alpha_{(i
l), t_{n_j}} g_j(u_{t_1},\dots,u_{t_{n_j}}) = 0.\]
Hence $f(a_1,\dots,a_n)=0,$ and $f$ is exact for
$A.$ \hfill $\Box$

\begin{lemma} \label{Kmu}
Let $\Gamma$ be a proper $\mathbb{Z}_2$T-ideal, and $A$ be a $\mathbb{Z}_2$PI-reduced algebra
such that $\mathrm{ind}_{\mathbb{Z}_2}(\Gamma)=\mathrm{ind}_{\mathbb{Z}_2}(A).$
Suppose that a graded polynomial $f$ satisfies the conditions
$f \notin \mathrm{Id}^{\mathbb{Z}_2}(A),$ and $f \in
K_{\hat{\mu}}(\Gamma)+\mathrm{Id}^{\mathbb{Z}_2}(A)$ for some $\hat{\mu} \in \mathbb{N}_0.$
Then $\mathbb{Z}_2 T[f] \bigcap S_{\mu}(A) \ne \emptyset$ holds
for any $\mu \in \mathbb{N}_0.$
\end{lemma}
\noindent {\bf Proof.} The full linearization $\tilde{f}$ of some
multihomogeneous component of $f$ also satisfies
$\tilde{f} \in K_{\hat{\mu}}(\Gamma)+\mathrm{Id}^{\mathbb{Z}_2}(A),$ \ $\tilde{f}
\notin \mathrm{Id}^{\mathbb{Z}_2}(A).$ Then by Lemma \ref{Exact4}
$\tilde{f}$ is exact for $A.$ And by Lemma \ref{gamma} we obtain that
$\emptyset \ne (\mathbb{Z}_2 T[\tilde{f}] \bigcap S_{\mu}(A)) \subseteq (\mathbb{Z}_2 T[f]
\bigcap S_{\mu}(A))$ for any $\mu \in \mathbb{N}_0.$ \hfill $\Box$

\section{Representable graded algebras.}

Let $R$ be a commutative associative unitary $F$-algebra.
Suppose that a $\mathbb{Z}/2 \mathbb{Z}$-graded $F$-algebra $A=A_{\bar{0}} \oplus A_{\bar{1}}$
has a structure of $R$-algebra satisfying $R A_\theta
\subseteq A_\theta,$ \ $\forall \theta \in \mathbb{Z}/2 \mathbb{Z}.$

\begin{definition}
Any $R$-linear mapping $\mathrm{tr}:A_{\bar{0}} \rightarrow R$ is
called trace on the $\mathbb{Z}/2 \mathbb{Z}$-graded $R$-algebra $A.$
\end{definition}
Observe that a trace on $A$ is not necessary symmetric (not necessary satisfies
$tr(ab)=tr(ba)$).

Denote by $\mathcal{S}$ the free associative commutative unitary algebra
generated by all symbols $\mathrm{tr}(u)$ for nonempty
associative noncommutative even monomials $u \in (F\langle
X^{\mathbb{Z}_2}\rangle)_{\bar{0}}$ over $X^{\mathbb{Z}_2}.$  We say that $FS\langle
X^{\mathbb{Z}_2}\rangle = F\langle X^{\mathbb{Z}_2}\rangle \otimes_F \mathcal{S}$ is free
$\mathbb{Z}/2 \mathbb{Z}$-graded algebra with trace. We assume that $(f \otimes s_1)
s_2 = s_2 (f \otimes s_1) = f \otimes (s_1 s_2)$ for all $f \in
F\langle X\rangle,$ \ $s_1, s_2 \in \mathcal{S}.$ The $\mathbb{Z}/2 \mathbb{Z}$-grading
on $FS\langle X^{\mathbb{Z}_2}\rangle$ is induced from $F\langle X^{\mathbb{Z}_2}\rangle$
($\mathcal{S}$ is supposed to be trivially graded). Elements of
$FS\langle X^{\mathbb{Z}_2}\rangle$ are called graded polynomials with
trace. Elements of $\mathcal{S}$ are called pure trace
polynomials.

We identify $F\langle X^{\mathbb{Z}_2}\rangle \otimes_F 1$ with
$F\langle X^{\mathbb{Z}_2}\rangle.$ The symbol $\otimes$ usually is omitted
in the notation of graded polynomials with trace. The
concept of degrees (homogeneity in the sense of
degree, multilinearity, alternating, etc.) of graded polynomials
with trace or pure trace polynomials is defined in similar way
to ordinary graded polynomials assuming $\deg_x
\mathrm{tr}(u)=\deg_x u$ for any $x \in X^{\mathbb{Z}_2}.$

The $\mathcal{S}$-linear function of trace
$\mathrm{tr}:(FS\langle X^{\mathbb{Z}_2}\rangle)_e \rightarrow \mathcal{S}$ is
defined on the $\mathcal{S}$-algebra $FS\langle X^{\mathbb{Z}_2}\rangle$ by the
formula \ $\mathrm{tr}(\sum_{i} u_i s_i)=\sum_{i} \mathrm{tr}(u_i)
s_i,$ where $u_i \in (F\langle X^{\mathbb{Z}_2}\rangle)_e$ are monomials on
$X^{\mathbb{Z}_2}$ of the even graded degree, $s_i \in \mathcal{S}.$

Let $A$ be a $\mathbb{Z}/2 \mathbb{Z}$-graded $R$-algebra with trace,
$f(x_{1},\dots,x_{n}) \in FS\langle X^{\mathbb{Z}_2}\rangle$ be a graded polynomial with trace.
$A$ satisfies the graded
identity with trace $f=0$ if $f(a_{1},\dots,a_{n})=0$ holds in $A$ for any
$a_{1}, \dots ,a_{n} \in A.$
The ideal of graded identities with trace of $A$
$\mathrm{SId}^{\mathbb{Z}_2}(A)=\{
f \in FS\langle X^{\mathbb{Z}_2}\rangle | \  A $ satisfies $ f=0 \}$ is
$\mathbb{Z}/2 \mathbb{Z}$-graded $\mathcal{S}$-ideal of $FS\langle X^{\mathbb{Z}_2}\rangle$ closed
under all graded endomorphisms of the algebra $FS\langle X^{\mathbb{Z}_2}\rangle,$ and
preserving the trace. $\mathrm{SId}^{\mathbb{Z}_2}(A)$ also satisfies the
condition $g \cdot \mathrm{tr}(f) \in \mathrm{SId}^{\mathbb{Z}_2}(A),$ \
for any $g \in FS\langle X^{\mathbb{Z}_2}\rangle,$ and for any $f \in
(\mathrm{SId}^{\mathbb{Z}_2}(A))_{\bar{0}}.$ Ideals of $FS\langle X^{\mathbb{Z}_2}\rangle$ with
these properties are called $\mathbb{Z}_2$TS-ideals.
Given a $\mathbb{Z}_2$TS-ideal $\widetilde{\Gamma},$ and polynomials with trace
$f, g \in FS\langle X^{\mathbb{Z}_2}\rangle$ we write
$f=g \ (\mathrm{mod} \ \widetilde{\Gamma})$ if $f-g \in
\widetilde{\Gamma}.$ $\mathbb{Z}_2 TS[\mathcal{V}]$ is the $\mathbb{Z}_2$TS-ideal generated by a set
$\mathcal{V} \subseteq FS\langle X^{\mathbb{Z}_2}\rangle.$

Let $\widetilde{\Gamma} \unlhd FS\langle
X^{\mathbb{Z}_2}\rangle$ be a $\mathbb{Z}_2$TS-ideal of
$FS\langle X^{\mathbb{Z}_2}\rangle$. Denote by
$I=\mathrm{Span}_F\{
\mathrm{tr}(f) v | f \in \widetilde{\Gamma}_{\bar{0}}, v \in \mathcal{S}
\} \unlhd \mathcal{S}$ the ideal of $\mathcal{S}$ generated by
all elements of the form $\mathrm{tr}(f)$ for all polynomials
with trace $f \in \widetilde{\Gamma}$ of the even graded degree.
Let $\bar{\mathcal{S}}=\mathcal{S}/I$ be the quotient algebra.
Then the quotient algebra $\overline{FS}\langle
X^{\mathbb{Z}_2}\rangle=FS\langle X^{\mathbb{Z}_2}\rangle /\widetilde{\Gamma}$ is $\mathbb{Z}/2 \mathbb{Z}$-graded and has the well-defined structure of $\bar{\mathcal{S}}$-algebra.
$\bar{\mathcal{S}}$-linear function $\mathrm{tr}:(FS\langle
X^{\mathbb{Z}_2}\rangle / \widetilde{\Gamma})_{\bar{0}} \rightarrow \bar{\mathcal{S}}$
is naturally defined by the equalities
$\mathrm{tr}(a+\widetilde{\Gamma})=\mathrm{tr}(a) + I$  for any
$a \in (FS\langle X^{\mathbb{Z}_2}\rangle)_{\bar{0}}.$
The ideal of graded identities with trace of $\overline{FS}\langle
X^{\mathbb{Z}_2}\rangle$ coincides with $\widetilde{\Gamma}.$ Moreover
$FS\langle X^{\mathbb{Z}_2}\rangle /\widetilde{\Gamma}$ is the
relatively free $\mathbb{Z}/2 \mathbb{Z}$-graded algebra with trace for the
$\mathbb{Z}_2$TS-ideal $\widetilde{\Gamma}.$

We also can consider the free associative $\mathbb{Z}/2 \mathbb{Z}$-graded algebra with
trace $FS\langle X_{\nu}^{\mathbb{Z}_2}\rangle$ and the relatively free $\mathbb{Z}/2 \mathbb{Z}$-graded
algebras with trace $FS\langle
X_{\nu}^{\mathbb{Z}_2}\rangle/(\widetilde{\Gamma} \bigcap FS\langle
X_{\nu}^{\mathbb{Z}_2}\rangle)$ of a finite rank $\nu \in \mathbb{N}.$

Let us take a finite dimensional $F$-superalgebra $A=B \oplus J$
with the Jacobson radical $J=J(A),$ and the semisimple part $B$. Then
the even component of the semisimple part
$B_{\bar{0}} = B \cap A_{\bar{0}}$ is a finite dimensional subalgebra of $A$
with $\dim_F B_{\bar{0}}=t_{\bar{0}},$ where $\mathrm{dims}_{\mathbb{Z}_2} A=(t_{\bar{0}},t_{\bar{1}}).$
Since the left regular representation $\mathfrak{T}:B_{\bar{0}}
\rightarrow M_{t_1}(F)$ of $B_{\bar{0}}$ is injective then $B_{\bar{0}}$ is
isomorphic to a subalgebra of the matrix algebra $M_{t_{\bar{0}}}(F).$

Therefore the trace $\mathrm{tr}:A_{\bar{0}} \rightarrow F$ on $A$ is well defined by
the rule
\begin{eqnarray} \label{Atrace}
\mathrm{tr}(a)=\mathrm{tr}(b+r)=\mathrm{Tr}(\mathfrak{T}(b)),
\quad a \in A_{\bar{0}}, \ b \in B_{\bar{0}}, \  r \in J_{\bar{0}},
\end{eqnarray}
where $\mathrm{Tr}$ is the usual trace of a linear
operator.

\begin{lemma} \label{Traceid1}
A finite dimensional $F$-superalgebra $A$ with the trace
(\ref{Atrace}) satisfies the $\mathbb{Z}/2 \mathbb{Z}$-graded identity with trace \
$\mathrm{tr}(z)f=\sum_{i=1}^{t_{\bar{0}}} f|_{x_i:=z x_i}.$ Where
$f(x_1,\dots,x_{t_1},Y) \in F\langle X^{\mathbb{Z}_2}\rangle$ is any
graded polynomial (without trace) of the
type $(\mathrm{dims}_{\mathbb{Z}_2} A,\mathrm{nd}(A)-1,1),$
alternating in the set $\{x_1,\dots,x_{t_{\bar{0}}} \} \subseteq X_{\bar{0}}.$ Here
$t_{\bar{0}}=\dim B_{\bar{0}},$ $z \in X_{\bar{0}}.$
\end{lemma}
\noindent {\bf Proof.} A polynomial $f$ of the type
$(\mathrm{dims}_{\mathbb{Z}_2} A,$ $\mathrm{nd}(A)-1,1)$ is $(\mathrm{dims}_{\mathbb{Z}_2}
A)$-alternating in at least one set of graded variables. Suppose that
$\{x_1, \dots, x_{t_{\bar{0}}} \} \subseteq X_{\bar{0}}$ is the part of the even
graded degree of this set. An elementary evaluation of
variables of the set $\{z, x_1, \dots, x_{t_{\bar{0}}} \}$ of the polynomial
$g=\mathrm{tr}(z)f-\sum_{i=1}^{t_{\bar{0}}} f|_{x_i:=z x_i}$ must be
semisimple. Moreover the set $\{x_1, \dots, x_{t_{\bar{0}}} \}$ must be
exchanged by pairwise different semisimple elements of the basis $D$
(\ref{basisD}) of degree $\bar{0}$. Otherwise, we get zero.

Suppose that $z=b \in B_{\bar{0}}$ in our evaluation.
Consider a polynomial $\tilde{f}(x_1, \dots, x_{t_{\bar{0}}},\dots)$ that is alternating
in a homogeneous set of variables of even degree $\{x_1, \dots, x_{t_{\bar{0}}} \} \subseteq X_{\bar{0}}.$
It can be directly checked (see also [\cite{BelRow}, Theorem J]) that for
an arbitrary linear operator $\mathfrak{K}: B_{\bar{0}} \rightarrow B_{\bar{0}}$
and for all pairwise distinct basic elements $b_1,\dots,b_{t_{\bar{0}}}$ of $B_{\bar{0}}$
the equality $Tr(\mathfrak{K}) \tilde{f}(b_1,\dots,b_{t_{\bar{0}}},\dots) =
\tilde{f}(\mathfrak{K}(b_1),\dots,b_{t_{\bar{0}}},\dots)+\dots+
\tilde{f}(b_1,\dots,\mathfrak{K}(b_{t_{\bar{0}}}),\dots)$
holds in $A$, the replacement of other variables being arbitrary.
Applying this observation to the linear operator $\mathfrak{T}(b)$ of the left multiplication
by the element $b$ we complete the proof.
\hfill $\Box$

Observe that in Lemma \ref{Traceid1} it is enough to consider semisimple
or radical evaluations of variables (not necessary elementary ones).

\begin{lemma} \label{Traceid2}
Let us take any $\beta \in \mathbb{N}_0^2,$ \ $\gamma \in
\mathbb{N}.$ Consider a graded polynomial without trace
$f(y_1,\dots,y_k) \in F\langle X^{\mathbb{Z}_2}\rangle$
of the type $(\beta;\gamma-1;1),$ and a pure trace polynomial $s(z_1,\dots,z_d) \in \mathcal{S}$ ($\{z_1,\dots,z_d\} \subseteq X^{\mathbb{Z}_2}$.) Then
there exists a graded polynomial without trace
$g_s(y_1,\dots,y_k,z_1,\dots,z_d) \in F\langle X^{\mathbb{Z}_2}\rangle$
such that $g_s \in \mathbb{Z}_2 T[f],$ and any finite
dimensional superalgebra $A$ with the parameters
$\mathrm{par}_{\mathbb{Z}_2}(A)=(\beta;\gamma)$ satisfies the graded identity
with trace
\begin{equation*}
s(z_1,\dots,z_d) \cdot f(y_1,\dots,y_k) -
g_s(y_1,\dots,y_k,z_1,\dots,z_d) =0.
\end{equation*}
\end{lemma}
\noindent {\bf Proof.} If
$\mathrm{par}_{\mathbb{Z}_2}(A)=(\beta;\gamma)$ then by Lemma \ref{Traceid1}
$A$ satisfies the $\mathbb{Z}/2 \mathbb{Z}$-graded
trace identity $\mathrm{tr}(z)f(w_1 x_1,w_2 x_2,\dots,w_{t_{\bar{0}}}
x_{t_{\bar{0}}},Y)-( f(z w_1 x_1, w_2 x_2,\dots,w_{t_{\bar{0}}} x_{t_{\bar{0}}},Y)+
\dots+f(w_1 x_1,w_2 x_2,\dots, z w_{t_{\bar{0}}} x_{t_{\bar{0}}},Y))=0.$
Where $w_i \in (F\langle X^{\mathbb{Z}_2}\rangle^{\#})_{\bar{0}}$
are arbitrary (possibly empty) graded monomials ($i=1,\dots,t_{\bar{0}}$),
and $f(x_1,\dots,x_{t_{\bar{0}}},Y) \in F\langle X^{\mathbb{Z}_2}\rangle$
is any graded polynomial of the type
$(\beta;\gamma-1;1),$ alternating in $\{x_1,\dots,x_{t_{\bar{0}}}\}
\subseteq X_{\bar{0}},$ \ $z \in X_{\bar{0}}.$

Let us take any graded monomials $u_i \in (F\langle X^{\mathbb{Z}_2}\rangle)_{\bar{0}}$
($i=1,\dots,n$), and any (possibly empty) graded monomials $v_j \in (F\langle
X^{\mathbb{Z}_2}\rangle^{\#})_{\bar{0}}$ ($j=1,\dots,t_{\bar{0}}$).
Then by induction on the number $n \in \mathbb{N}_0$
of monomials $u_i$  there exist a natural
$\widetilde{n},$ and graded monomials $\widetilde{v}_{l j} \in
(F\langle X^{\mathbb{Z}_2}\rangle^{\#})_{\bar{0}}$ (possibly empty) such that
$\mathrm{tr}(u_1) \cdots \mathrm{tr}(u_n) f(v_1
x_1,\dots,v_{t_{\bar{0}}} x_{t_{\bar{0}}},Y)=$ \\ $\sum_{l=1}^{\widetilde{n}}
f(\widetilde{v}_{l 1} x_1,\dots,\widetilde{v}_{l t_{\bar{0}}} x_{t_{\bar{0}}},Y)$ \
$(\mathrm{mod} \ \mathrm{SId}^{\mathbb{Z}_2}(A)).$ Moreover the monomials
$\widetilde{v}_{l j}$
depend on the same variables as the monomial  $u_1 \cdots u_n v_1 \cdots v_{t_{\bar{0}}}$.

Hence for any pure trace polynomial $s(z_1,\dots,z_d)=\sum_{(j)} \alpha_{(j)}
\mathrm{tr}(u_{j_1}) \cdots \mathrm{tr}(u_{j_n}) \in \mathcal{S}$
the algebra
$A$ satisfies the identity $$s(z_1,\dots,z_d) \cdot
f(x_1,\dots,x_{t_{\bar{0}}},Y) - g_s(x_1,\dots,x_{t_{\bar{0}}},Y,z_1,\dots,z_d)=0.$$ Where
$g_s \in \mathbb{Z}_2 T[f]$ is some graded polynomial that does not depend on $A.$
Here $u_j \in (F\langle
z_1,\dots,z_d \rangle)_{\bar{0}}$ are monomials, $\alpha_{(j)} \in F.$
\hfill $\Box$

Let $A$ be a finite dimensional $F$-superalgebra with the
semisimple part $B=B_{\bar{0}} \oplus B_{\bar{1}},$ and the
Jacobson radical $J=J_{\bar{0}} \oplus J_{\bar{1}}$.
Let us denote $\dim B_\theta = t_\theta,$ $\dim J_\theta = q_\theta$ for any
$\theta \in \mathbb{Z}/2 \mathbb{Z}$.

Given a number $\nu \in \mathbb{N}$ take a set
$\Lambda_{\nu}=\{ \lambda_{\theta i j} | \theta \in \mathbb{Z}/2 \mathbb{Z}, 1 \le i \le
\nu, 1 \le j \le t_\theta+q_\theta \}.$
Consider the free commutative associative unitary algebra
$F[\Lambda_{\nu}]^{\#}$ generated by
$\Lambda_{\nu}$, and the associative algebra
$\mathcal{P}_{\nu}(A)=F[\Lambda_{\nu}]^{\#} \otimes_F A.$
The algebra $\mathcal{P}_{\nu}(A)$ has the structure of $F[\Lambda_{\nu}]
^{\#}$-module defined by \ $a \cdot f = f \cdot a  = f \cdot (\sum_{i} f_i
\otimes a_i)=\sum_{i} (f f_i) \otimes a_i,$ \ for any $f, f_i \in
F[\Lambda_{\nu}]^{\#},$ \ $a_i \in A,$ \ $a=\sum_{i} f_i \otimes
a_i \in \mathcal{P}_{\nu}(A).$

The $\mathbb{Z}/2 \mathbb{Z}$-grading on $\mathcal{P}_{\nu}(A)$ is induced from $A$
assuming that $F[\Lambda_{\nu}]^{\#}$ is trivially graded.
Define an $F[\Lambda_{\nu}]^{\#}$-linear map
$\mathrm{tr}:(\mathcal{P}_{\nu}(A))_{\bar{0}} \rightarrow
F[\Lambda_{\nu}]^{\#}$ by the equalities
\begin{eqnarray} \label{Ptrace}
\mathrm{tr}(a)=\mathrm{tr}(\sum_{i=1}^{t_{\bar{0}}} f_i \otimes b_{i {\bar{0}}} +
\sum_{j=1}^{q_{\bar{0}}} \tilde{f}_j \otimes r_{j {\bar{0}}})=\sum_{i=1}^{t_{\bar{0}}} f_i
\ \mathrm{Tr}(\mathfrak{T}(b_{i {\bar{0}}})).
\end{eqnarray}
$\mathrm{tr}$ is well defined in $\mathcal{P}_{\nu}(A).$ Here $b_{i {\bar{0}}} \in B_{\bar{0}},$
\ $r_{j e} \in (J(A))_{\bar{0}},$ \ $f_i, \tilde{f}_j \in
F[\Lambda_{\nu}]^{\#}.$ \  $\mathfrak{T}$ is the left regular
representation of $B_{\bar{0}},$ and $\mathrm{Tr}$ is the usual trace.

\begin{lemma} \label{HamKel}
$\mathcal{P}_{\nu}(A)$ satisfies the graded identity with trace
$(\mathcal{X}_{t_{\bar{0}}}(x) \cdot x)^{\mathrm{nd}(A)}=0.$ Here
$\mathcal{X}_{t_{\bar{0}}}(x)$ is the Cayley-Hamilton polynomial of the
degree $t_{\bar{0}},$ \ $x \in X_{\bar{0}}.$
\end{lemma}
\noindent {\bf Proof.}
The algebra $B_{\bar{0}}$ with the trace (\ref{Atrace}) satisfies
the identity with trace $\mathcal{X}_{t_{\bar{0}}}(x) \cdot x=0$
(\cite{Proc1}, \cite{Rasm1}).
$F[\Lambda_{\nu}]^{\#}$ is commutative
non-nilpotent algebra. $F[\Lambda_{\nu}]^{\#} \otimes_F (J(A))_{\bar{0}}$
is a nilpotent ideal of
$(\mathcal{P}_{\nu}(A))_{\bar{0}}$ of the degree $\mathrm{nd}(A).$
Hence the neutral component
$(\mathcal{P}_{\nu}(A))_{\bar{0}}=(F[\Lambda_{\nu}]^{\#} \otimes_F B_{\bar{0}})
\oplus (F[\Lambda_{\nu}]^{\#} \otimes_F (J(A))_{\bar{0}})$ satisfies the
full linearization of the identity $(\mathcal{X}_{t_{\bar{0}}}(x) \cdot
x)^{\mathrm{nd}(A)}=0.$ \hfill
$\Box$

Let  $\{\hat{b}_{\theta 1}, \dots,
\hat{b}_{\theta t_\theta}\}$ be a basis of the graded component
$B_\theta$ ($\theta \in \mathbb{Z}/2 \mathbb{Z}$) of the semisimple part $B$ of $A,$ \ $\dim B_\theta =
t_\theta.$ Also let
$\{\hat{r}_{\theta 1}, \dots, \hat{r}_{\theta
q_\theta}\}$ be a basis of the graded component $J_\theta$ of the
Jacobson radical $J=J(A),$ \ $\dim J_\theta = q_\theta.$ Recall
that for a superalgebra with elementary decomposition all these bases
may be chosen in the set $D \bigcup U$
((\ref{basisD}), (\ref{basisU}), Lemma \ref{Pierce}). Take
elements
\begin{eqnarray} \label{genset}
y_{\theta i}=\sum_{j=1}^{t_\theta} \lambda_{\theta i j} \otimes
\hat{b}_{\theta j} + \sum_{j=1}^{q_\theta} \lambda_{\theta i
j+t_\theta} \otimes \hat{r}_{\theta j} \  \in  \mathcal{P}_{\nu}(A),
\qquad \theta \in \mathbb{Z}/2 \mathbb{Z}, \  \ 1 \le i \le \nu.
\end{eqnarray}
The elements $y_{\theta i}$ are $\mathbb{Z}/2 \mathbb{Z}$-homogeneous of degree $\theta$
($\theta \in \mathbb{Z}/2 \mathbb{Z}; \ 1 \le i \le \nu$). Given a positive integer  $\nu$
consider
the $F$-subalgebra $\mathcal{F}_{\nu}(A)=\langle y_{\theta i} | \theta
\in \mathbb{Z}/2 \mathbb{Z}, \ 1 \le i \le \nu \rangle$ of $\mathcal{P}_{\nu}(A)$
generated by the set $Y^{\mathbb{Z}_2}_{\nu}=\{ y_{\theta i} | \theta \in \mathbb{Z}/2 \mathbb{Z}, \
1 \le i \le \nu \}.$  $\mathcal{F}_{\nu}(A)$ is
$\mathbb{Z}/2 \mathbb{Z}$-graded.

Any map $\varphi$ of generators to
arbitrary homogeneous elements $\widetilde{a}_{\theta i} \in A_\theta$ \
($\widetilde{\alpha}_{\theta i j} \in F$)
\begin{equation} \label{hom1}
\varphi:y_{\theta i} \mapsto \widetilde{a}_{\theta
i}=\sum_{j=1}^{t_\theta} \widetilde{\alpha}_{\theta i j}
\hat{b}_{\theta j} + \sum_{j=1}^{q_\theta}
\widetilde{\alpha}_{\theta i j+t_\theta} \hat{r}_{\theta j} \in
A_\theta \qquad (\theta \in \mathbb{Z}/2 \mathbb{Z}, \ \  i=1,\dots,\nu)
\end{equation}
can be extended to the graded homomorphism of $F$-superalgebras
$\varphi:\mathcal{F}_{\nu}(A) \rightarrow A$, also inducing the
graded homomorphism $\widetilde{\varphi}:\mathcal{P}_{\nu}(A)
\rightarrow A$ defined by the following equalities
\begin{eqnarray} \label{hom2}
\widetilde{\varphi}((\lambda_{\theta_1 i_1 j_1} \cdots
\lambda_{\theta_k i_k j_k})\otimes a)=(\widetilde{\alpha}_{\theta_1
i_1 j_1} \cdots \widetilde{\alpha}_{\theta_k i_k j_k}) \cdot a
\qquad \forall a \in A.
\end{eqnarray}
The homomorphism $\widetilde{\varphi}$ preserves the trace
defined by (\ref{Ptrace}) on $\mathcal{P}_{\nu}(A)$ and
by (\ref{Atrace}) on $A.$

Elements of $\mathcal{F}_{\nu}(A)$
are called quasi-polynomials in the variables $Y^{\mathbb{Z}_2}_{\nu}.$ Products
of the generators $y_{\theta i} \in Y^{\mathbb{Z}_2}_{\nu}$ of the algebra
$\mathcal{F}_{\nu}(A)$ are called quasi-monomials. We have
also $\mathrm{Id}^{\mathbb{Z}_2}(\mathcal{F}_{\nu}(A)) \supseteq
\mathrm{Id}^{\mathbb{Z}_2}(\mathcal{P}_{\nu}(A))=\mathrm{Id}^{\mathbb{Z}_2}(A)$ for any $\nu
\in \mathbb{N}.$

$\mathcal{F}_{\nu}(A)$ is a finitely generated PI-algebra. By
Shirshov's height theorem  \cite{Shirsh} $\mathcal{F}_{\nu}(A)$
has a finite height and a finite Shirshov's basis.
Shirshov's basis of an algebra always can be chosen in the set of
monomials over the generators (\cite{BelRow}, \cite{Shirsh}). Thus
we can suppose that a Shirshov's basis of $\mathcal{F}_{\nu}(A)$
consists of homogeneous in the grading elements of
$\mathcal{F}_{\nu}(A).$ More precisely, there exist a positive integer
$\mathcal{H},$ and homogeneous in the grading
elements $w_1,\dots,w_d \in \mathcal{F}_{\nu}(A)$ such that
any element $u \in \mathcal{F}_{\nu}(A)$ has the form
$u=\sum_{(i)=(i_1,\dots,i_k)}
 \alpha_{(i)} \  w_{i_1}^{c_1} \dots w_{i_k}^{c_k},$
where $k \le \mathcal{H},$ \  $\{ i_1,\dots,i_k \} \subseteq
\{1,\dots,d\},$ \  $c_j \in \mathbb{N},$ \  $\alpha_{(i)} \in F.$

Consider the polynomials
$\hat{\mathfrak{s}}_{i l}=\mathrm{tr}(w_i^{2 l}) \in
F[\Lambda_{\nu}]^{\#}$ ($i=1,\dots,d,$ \ $l=1,\dots,t_{\bar{0}}$).
Then
$\widehat{F}=F[\hat{\mathfrak{s}}_{i l} \ | \ 1 \le i \le d; \ 1
\le l \le t_{\bar{0}} \ ]^{\#}$ is
the associative commutative non-graded $F$-subalgebra
of $F[\Lambda_{\nu}]^{\#}$ with the unit
generated by $\{ \hat{\mathfrak{s}}_{i l} \}$, and by
the unit of $F[\Lambda_{\nu}]^{\#}.$

Take the $\mathbb{Z}/2 \mathbb{Z}$-graded $\widehat{F}$-subalgebra
$\mathcal{T}_{\nu}(A)=\widehat{F} \mathcal{F}_{\nu}(A)$ of
$\mathcal{P}_{\nu}(A).$ Then $\mathcal{F}_{\nu}(A)$ is a graded
subalgebra of $\mathcal{T}_{\nu}(A).$ An arbitrary
map of type (\ref{hom1}) can be uniquely extended to a
graded homomorphism from $\mathcal{T}_{\nu}(A)$ to $A$ preserving the
traces (it is the restriction of
$\widetilde{\varphi}$ defined by (\ref{hom2}) onto $\mathcal{T}_{\nu}(A)$ ).

Since for any $i=1,\dots,d$ we have $w_i^2 \in (\mathcal{P}_{\nu}(A))_{\bar{0}}$ \
then it follows from Lemma \ref{HamKel} that all elements $w_i^2$ are
algebraic of degree $\mathrm{nd}(A)(t_{\bar{0}}+1)$ over
$\widehat{F}.$ Therefore, by Shirshov's height theorem,
$\mathcal{T}_{\nu}(A)$ is finitely generated
$\widehat{F}$-module, where $\widehat{F}$ is Noetherian. By
theorem of Beidar \cite{Beid} the algebra $\mathcal{T}_{\nu}(A)$
is representable.

Let $V \subseteq F\langle X^{\mathbb{Z}_2} \rangle$ be a set of graded
polynomials. We denote by $V(\mathcal{T}_{\nu}(A)) \unlhd
\mathcal{T}_{\nu}(A)$ the verbal ideal generated by results of all
appropriate substitutions of $\mathbb{Z}/2 \mathbb{Z}$-homogeneous elements of
$\mathcal{T}_{\nu}(A)$ to any graded polynomial from $V.$

\begin{lemma} \label{VerbId}
Given a set $V \subseteq F\langle X^{\mathbb{Z}_2} \rangle,$ and any
$\nu \in \mathbb{N}$ the verbal ideal
\begin{eqnarray*}
&&V(\mathcal{T}_{\nu}(A))=\mathrm{Span}_F\{ \ \mathfrak{s} \cdot
v_1 \tilde{f}(u_1,\dots,u_n) v_2 \ | \
\mathfrak{s} \in \widehat{F}; \quad u_i \in \mathcal{F}_{\nu}(A)
\mbox{ are } \nonumber\\ &&\mbox{quasi-monomials; } \quad
v_1, v_2 \in \mathcal{F}_{\nu}(A)^{\#} \
\mbox{ are quasi-monomials, possibly empty; } \nonumber \\ &&\tilde{f} \
\mbox{ is the full linearization of a multihomogeneous component of any  } f \in V
\}
\end{eqnarray*}
is a graded $\widehat{F}$-closed ideal of $\mathcal{T}_{\nu}(A).$
The
quotient algebra
$\overline{\mathcal{T}}_{\nu}(A,V)=\mathcal{T}_{\nu}(A)/V(\mathcal{T}_{\nu}(A))$
is a representable $\widehat{F}$-superalgebra.
The ideal of graded identities of $\overline{\mathcal{T}}_{\nu}(A,V)$ satisfies
\ $\mathrm{Id}^{\mathbb{Z}_2}(\overline{\mathcal{T}}_{\nu}(A,V)) \supseteq
\widetilde{\Gamma}+\mathrm{Id}^{\mathbb{Z}_2} \bigl(F\langle X_{\nu}^{\mathbb{Z}_2}
\rangle/(\widetilde{\Gamma} \cap F\langle X_{\nu}^{\mathbb{Z}_2}
\rangle ) \bigr),$ where $\widetilde{\Gamma}=\mathrm{Id}^{\mathbb{Z}_2}(A)+\mathbb{Z}_2 T[V].$
\end{lemma}
\noindent {\bf Proof.} It is clear that in case of
characteristic zero any verbal ideal of an algebra can
be generated by the results of all appropriate substitutions to the
full linearizations of multihomogeneous components
of polynomials from the given set. Particularly,
$V(\mathcal{T}_{\nu}(A))$ is generated by $\mathbb{Z}/2 \mathbb{Z}$-homogeneous elements,
and hence, this is a graded ideal.

Any homogeneous element of $\mathcal{T}_{\nu}(A)$
has the form $c=\sum_{i} \mathfrak{s}_i \cdot u_i,$ where
$\mathfrak{s}_i \in \widehat{F},$ \ $u_i \in \mathcal{F}_{\nu}(A)$
are quasi-monomials, $\deg_{\mathbb{Z}/2 \mathbb{Z}} u_i=\deg_{\mathbb{Z}/2 \mathbb{Z}} c$
for all $i.$ If
$\tilde{f}$ is a multilinear graded polynomial then
$\tilde{f}(c_1,\dots,c_n)=\tilde{f}(\sum_{i_1} \mathfrak{s}_{1 i_1}
\cdot u_{1 i_1},\dots,$ $\sum_{i_n} \mathfrak{s}_{n i_n} \cdot u_{n
i_n})=\sum_{(i_1,\dots,i_n)} \bigl( \prod_{l=1}^{n} \mathfrak{s}_{l
i_l} \bigr)
 \tilde{f}(u_{1 i_1},\dots,u_{n i_n})$
for any homogeneous elements $c_1, \dots, c_n \in
\mathcal{T}_{\nu}(A).$ Here $u_{l i_l} \in \mathcal{F}_{\nu}(A)$ are
quasi-monomials of appropriate graded degrees, $\mathfrak{s}_{l
i_l} \in \widehat{F}.$ We have also
$\mathfrak{s} \cdot \tilde{f}(c_1,\dots,c_n)=\tilde{f}(\mathfrak{s}
\cdot c_1,\dots,c_n)=\tilde{f}(c'_1,\dots,c_n) \in
V(\mathcal{T}_{\nu}(A))$ for any $\mathfrak{s} \in \widehat{F}.$

Then the quotient algebra $\overline{\mathcal{T}}_{\nu}(A,V)$ is
$\mathbb{Z}/2 \mathbb{Z}$-graded is an $\widehat{F}$-module.
$\overline{\mathcal{T}}_{\nu}(A,V)$ is
finitely generated over $\widehat{F},$ as well as
$\mathcal{T}_{\nu}(A),$ and is also representable.
$\mathcal{T}_{\nu}(A)$ is a graded $F$-subalgebra of
$\mathcal{P}_{\nu}(A).$ Therefore we have
$\mathrm{Id}^{\mathbb{Z}_2}(\overline{\mathcal{T}}_{\nu}(A,V)) \supseteq
\mathrm{Id}^{\mathbb{Z}_2}(\mathcal{T}_{\nu}(A)) \supseteq
\mathrm{Id}^{\mathbb{Z}_2}(\mathcal{P}_{\nu}(A)) = \mathrm{Id}^{\mathbb{Z}_2}(A).$ It is
clear also that $V \subseteq
\mathrm{Id}^{\mathbb{Z}_2}(\overline{\mathcal{T}}_{\nu}(A,V)),$ hence
$\widetilde{\Gamma}=\mathrm{Id}^{\mathbb{Z}_2}(A)+\mathbb{Z}_2T[V] \subseteq
\mathrm{Id}^{\mathbb{Z}_2}(\overline{\mathcal{T}}_{\nu}(A,V)).$

Let us denote the $\mathbb{Z}_2$T-ideals $\Gamma_1=\mathrm{Id}^{\mathbb{Z}_2} \bigl(F\langle X_{\nu}^{\mathbb{Z}_2}
\rangle/( \widetilde{\Gamma} \cap F\langle X_{\nu}^{\mathbb{Z}_2}
\rangle ) \bigr),$ and $\Gamma_2=\mathrm{Id}^{\mathbb{Z}_2}(\overline{\mathcal{T}}_{\nu}(A,V)).$
From the arguments above we have $\Gamma_1(F\langle X_{\nu}^{\mathbb{Z}_2}\rangle) \subseteq
\widetilde{\Gamma} \subseteq \Gamma_2.$
 Then for a multilinear graded polynomial $f(x_1,\dots,x_n) \in
\Gamma_1,$ for any homogeneous
elements $u_i \in \mathcal{F}_{\nu}(A)$ ($\deg_G u_i=\deg_G
x_i$), and for any $\mathfrak{s}_i \in \widehat{F}$
we obtain
\begin{eqnarray*}
f(\mathfrak{s}_1 u_1,\dots,\mathfrak{s}_n u_n)=\mathfrak{s}_1 \cdots
\mathfrak{s}_n f(u_1,\dots,u_n) \in \widehat{F}
\Gamma_1(\mathcal{F}_{\nu}(A)) \subseteq \widehat{F}
\Gamma_2(\mathcal{F}_{\nu}(A)) \subseteq V(\mathcal{T}_{\nu}(A)),
\end{eqnarray*}
since $\mathrm{grk}(\mathcal{F}_{\nu}(A))=\nu.$ Therefore $f \in
\Gamma_2,$ and $\Gamma_1 \subseteq \Gamma_2.$ \hfill $\Box$

\begin{definition} \label{Homog}
We say that a subset $V \subseteq F\langle X^{\mathbb{Z}_2} \rangle$ is
multihomogeneous if for any $f \in V$ \ $V$ contains all multihomogeneous components of $f.$
\end{definition}

\begin{lemma} \label{GS1}
Let $A=A_1 \times \dots \times A_\rho$
be the direct product of arbitrary finite dimensional
superalgebras $A_1, \dots, A_\rho.$ Suppose that
$V \subseteq F\langle X^{\mathbb{Z}_2} \rangle$ is a multihomogeneous set, and
$\nu$ is any positive integer. Let us take any $f(z_1,\dots,z_n)
\in \mathrm{Id}^{\mathbb{Z}_2}(\overline{\mathcal{T}}_{\nu}(A,V)),$ and any
homogeneous polynomials $h_1,\dots,h_n \in F\langle
X_{\nu}^{\mathbb{Z}_2} \rangle$ with $\deg_{\mathbb{Z}_2} h_l = \deg_{\mathbb{Z}_2} z_l.$
Then the equality
$f(h_1, \dots,h_n)=\sum_{j} \mathfrak{s}_j
\cdot v_{j 1} \tilde{f}_j(u_{j 1},\dots,u_{j n})
v_{j 2} \ (\mathrm{mod } \ \mathrm{SId}^{\mathbb{Z}_2}(A_i))$ holds for
any $i=1,\dots,\rho.$ Here $\tilde{f}_j$ are the full
linearizations of some polynomials $f_j \in V;$ \ $u_{j l},
v_{j l} \in F\langle X_{\nu}^{\mathbb{Z}_2} \rangle$ are
monomials, $v_{j l}$ may be empty; and $\mathfrak{s}_j
\in \mathcal{S}$ are pure trace graded polynomials in the
variables $X_{\nu}^{\mathbb{Z}_2}.$
\end{lemma}
\noindent {\bf Proof.}
Given a graded polynomial $f(z_1,\dots,z_n) \in
\mathrm{Id}^{\mathbb{Z}_2}(\overline{\mathcal{T}}_{\nu}(A,V)),$ and arbitrary
homogeneous polynomials $h_1,\dots,h_n \in F\langle X_{\nu}^{\mathbb{Z}_2}
\rangle$ of degrees according to variables of $f,$ we have
$f(\tilde{h}_1,\dots,\tilde{h}_n) \in V(\mathcal{T}_{\nu}(A)).$ Here the quasi-polynomial
$\tilde{h}_i=h_i(y_1,\dots,y_{2 \nu})$ is obtained by replacement of the variables
$x_l \in X_{\nu}^{\mathbb{Z}_2}$ by the corresponding elements
$y_l \in Y_{\nu}^{\mathbb{Z}_2}$ (\ref{genset}). Since $V$ is a
multihomogeneous set then by Lemma \ref{VerbId} in the algebra
$\mathcal{T}_{\nu}(A)$ we obtain the equality
$f(\tilde{h}_1,\dots,\tilde{h}_n)= \sum_{j}
\tilde{\mathfrak{s}}_j \cdot \tilde{v}_{j 1} \tilde{f}_j(\tilde{u}_{j 1},\dots,
\tilde{u}_{j n}) \tilde{v}_{j 2},$
where $\tilde{f}_j$ are the full linearizations of polynomials $f_j
\in V,$ \ $\tilde{u}_{j l}=u_{j l}(y_1,\dots,y_{m \nu}) \in \mathcal{F}_{\nu}(A)$
are quasi-monomials, \ $\tilde{\mathfrak{s}}_j=\mathfrak{s}_j(y_1,\dots,y_{m \nu}) \in
\widehat{F}$ are pure trace quasi-polynomials, $\tilde{v}_{j l}=v_{j l}(y_1,\dots,y_{m \nu}) \in \mathcal{F}_{\nu}(A)$ are also quasi-monomials, possibly empty.

An arbitrary graded map $\varphi:Y^{\mathbb{Z}_2}_{\nu} \rightarrow
A_i$ from the generating set (\ref{genset}) of
$\mathcal{F}_{\nu}(A)$ into any subalgebra $A_i \subseteq A$ ($i=1,\dots,\rho$)
can be extended to the graded homomorphism
$\widetilde{\varphi}:\mathcal{T}_{\nu}(A) \rightarrow A_i$ preserving
the trace. Thus the equality
\begin{eqnarray*}
&&f(h_1(a_1,\dots,a_{2 \nu}),\dots,h_n(a_1,\dots,a_{2
\nu}))-\sum_{j} \mathfrak{s}_j(a_1,\dots,a_{2 \nu}) \cdot
v_{j 1}(a_1,\dots,a_{2 \nu}) \times \\
&&\tilde{f}_j(u_{j 1}(a_1,\dots,a_{2 \nu}),\dots,u_{j
n}(a_1,\dots,a_{2 \nu})) \  v_{j 2}(a_1,\dots,a_{2
\nu})=0
\end{eqnarray*}
holds in $A_i$ for any elements \
$a_1,\dots,a_{2 \nu} \in A_i$ of appropriate graded
degrees. Therefore
\begin{eqnarray*}
&&f(h_1(x_1,\dots,x_{2 \nu}), \dots,h_n(x_1,\dots,x_{2
\nu}))-\sum_{j} \mathfrak{s}_j(x_1,\dots,x_{2 \nu}) \cdot
v_{j 1}(x_1,\dots,x_{2 \nu}) \times \\
&&\tilde{f}_j(u_{j 1}(x_1,\dots,x_{2 \nu}),\dots,u_{j
n}(x_1,\dots,x_{2 \nu})) \  v_{j 2}(x_1,\dots,x_{2
\nu})=0
\end{eqnarray*}
is a graded identity with trace of the algebra
$A_i,$ for any $i=1,\dots,\rho.$   \hfill $\Box$

\begin{lemma} \label{GS}
Suppose that $A_1,\dots,A_\rho$ are
any $\mathbb{Z}_2$PI-reduced algebras  with
$\mathrm{ind}_{\mathbb{Z}_2}(A_i)=\kappa$ for all $i=1,\dots,\rho.$ Given a
subset $V \subseteq \bigcup_{i=1}^{\rho} S_{\mu}(A_i)$ (for any
$\mu \ge 1$), and a positive integer $\nu,$ there exists an
$F$-finite dimensional superalgebra $C_\nu$ such
that $\mathrm{Id}^{\mathbb{Z}_2}(C_\nu) = \mathrm{Id}^{\mathbb{Z}_2} \Bigl(F\langle
X_{\nu}^{\mathbb{Z}_2} \rangle/\bigl((\bigcap_{i=1}^\rho
\mathrm{Id}^{\mathbb{Z}_2}(A_i)+\mathbb{Z}_2 T[V]) \bigcap F\langle X_{\nu}^{\mathbb{Z}_2} \rangle
\bigr) \Bigr).$
\end{lemma}
\noindent {\bf Proof.}
Let us take $A=A_1 \times \dots \times A_\rho,$ and any $\nu \in \mathbb{N}.$
By Lemma \ref{dirsumKP} we have $\bigcup_{i=1}^{\rho}
S_{\mu}(A_i)=S_{\mu}(A)$ for any $\mu,$ hence $V \subseteq S_{\mu}(A).$
Observe that a set of boundary polynomials is always multihomogeneous.
Then for any $f(z_1,\dots,z_n) \in
\mathrm{Id}^{\mathbb{Z}_2}(\overline{\mathcal{T}}_{\nu}(A,V)),$ and any
homogeneous polynomials $h_1,\dots,h_n \in F\langle X_{\nu}^{\mathbb{Z}_2}
\rangle$ of appropriate graded degrees by Lemma \ref{GS1}
we have
\[ f(h_1, \dots,h_n)=\sum_{j} \mathfrak{s}_j v_{j 1}
\tilde{f}_j(u_{j 1},\dots,u_{j n}) v_{j 2} \
(\mathrm{mod } \ \  \mathrm{SId}^{\mathbb{Z}_2}(A_i)) \qquad \forall \  i=1,\dots,\rho.
\]
Where $\tilde{f}_j \in \mathbb{Z}_2 T[V]$ are multilinear graded polynomials;
\  $u_{j l}, v_{j l} \in F\langle
X_{\nu}^{\mathbb{Z}_2} \rangle$ are monomials, $v_{j l}$ may be
empty; \  $\mathfrak{s}_j \in \mathcal{S},$ \  $\mathfrak{s}_j$ depends on $X_{\nu}^{\mathbb{Z}_2}.$

$A_i$ is a $\mathbb{Z}_2 $PI-reduced superalgebra, hence from Lemmas \ref{gamma}, \ref{dirsum}
it follows that
$\mathrm{par}_G(A_i)=\mathrm{ind}_G(A_i)=
\mathrm{ind}_G(A)=\kappa=(\beta;\gamma)$ (\ $\forall
i=1,\dots,\rho$).  Any polynomial $\tilde{f}_j$ has the type
$(\beta;\gamma-1;\mu)$ (as the multilinearization of the polynomial
$f_j \in V \subseteq S_\mu(A)$). Then by Lemma
\ref{Traceid2} there exists a graded traceless polynomial $\tilde{g}_j \in
\mathbb{Z}_2 T[\tilde{f}_j] \cap F\langle X_{\nu}^{\mathbb{Z}_2} \rangle \subseteq
\mathbb{Z}_2 T[V] \cap F\langle X_{\nu}^{\mathbb{Z}_2} \rangle$ such that
$\mathfrak{s}_j \ \tilde{f}_j(v_{j
1},\dots,v_{j n})\  = \
\tilde{g}_j \  (\mathrm{mod } \
\mathrm{SId}^{\mathbb{Z}_2} (A_i))$ \ \ for any $i=1,\dots,\rho.$ Therefore
\[ f(h_1, \dots,h_n)=\sum_{j} v_{j 1}
\tilde{g}_j v_{j 2} \
(\mathrm{mod } \ \mathrm{SId}^{\mathbb{Z}_2}(A_i))
\qquad \forall \  i=1,\dots,\rho,
\]
where $g=\sum_{j}
v_{j 1} \tilde{g}_j v_{j 2} \in \mathbb{Z}_2 T[V]
\bigcap F\langle X_{\nu}^{\mathbb{Z}_2} \rangle.$ Hence we obtain that $f(h_1,\dots,h_n)-g \in$
$\mathrm{SId}^{\mathbb{Z}_2}(A_i) \cap F\langle X_{\nu}^{\mathbb{Z}_2} \rangle =$
$\mathrm{Id}^{\mathbb{Z}_2}(A_i) \cap F\langle X_{\nu}^{\mathbb{Z}_2} \rangle$ for all
$i=1,\dots,\rho.$ Thus $f(h_1,\dots,h_n) \in \bigl(\mathrm{Id}^{\mathbb{Z}_2}(A)
+ \mathbb{Z}_2 T[V] \bigr) \bigcap F\langle X_{\nu}^{\mathbb{Z}_2} \rangle,$ \ and $f \in
\mathrm{Id}^{\mathbb{Z}_2} \Bigl(F\langle X_{\nu}^{\mathbb{Z}_2}
\rangle/\bigl((\mathrm{Id}^{\mathbb{Z}_2}(A)+\mathbb{Z}_2 T[V]) \bigcap F\langle X_{\nu}^{\mathbb{Z}_2}
\rangle \bigr) \Bigr).$

By Lemma \ref{VerbId} we have also that \  $\mathrm{Id}^{\mathbb{Z}_2} \Bigl(F\langle X_{\nu}^{\mathbb{Z}_2}
\rangle/\bigl((\mathrm{Id}^{\mathbb{Z}_2}(A)+\mathbb{Z}_2 T[V]) \bigcap F\langle X_{\nu}^{\mathbb{Z}_2}
\rangle \bigr) \Bigr) \subseteq$
$\mathrm{Id}^{\mathbb{Z}_2}(\overline{\mathcal{T}}_{\nu}(A,V)).$ The superalgebra
$\overline{\mathcal{T}}_{\nu}(A,V)$ is representable.
Hence by Lemma \ref{Repr} there exists an $F$-finite dimensional
superalgebra $C_\nu$ such that
\[\mathrm{Id}^{\mathbb{Z}_2}(C_\nu)=\mathrm{Id}^{\mathbb{Z}_2}(\overline{\mathcal{T}}_{\nu}(A,V))=\mathrm{Id}^{\mathbb{Z}_2}
\Bigl(F\langle X_{\nu}^{\mathbb{Z}_2} \rangle/\bigl((\mathrm{Id}^{\mathbb{Z}_2}(A)+\mathbb{Z}_2 T[V])
\bigcap F\langle X_{\nu}^{\mathbb{Z}_2} \rangle \bigr) \Bigr).\] \hfill $\Box$

\section{Graded identities of finitely generated PI-algebras.}

\begin{lemma} \label{Main}
Let $F$ be a field of characteristic 0. Let $\Gamma$  be a non-trivial ideal of graded identities of
a finitely generated associative PI-superalgebra over $F.$
Then there exists a finite dimensional associative
$F$-superalgebra $\widetilde{A}$ satisfying the conditions $\mathrm{Id}^{\mathbb{Z}_2}(\widetilde{A})
\subseteq
\Gamma,$ \ $\mathrm{ind}_{\mathbb{Z}_2}(\Gamma)=\mathrm{ind}_{\mathbb{Z}_2}(\widetilde{A}),$ and
$S_{\hat{\mu}}(\mathcal{O}(\widetilde{A})) \cap \Gamma =\emptyset$ for some
$\hat{\mu} \in \mathbb{N}_0.$
\end{lemma}
\noindent {\bf Proof.} $\Gamma$ contains a non-trivial T-ideal $\widetilde{\Gamma}$
of ordinary non-graded identities of a finitely generated associative PI-algebra.
By \cite{Lewin} $\widetilde{\Gamma}$ also contains the T-ideal of some
finite dimensional non-graded algebra $C^{\mathfrak{od}}.$
It is clear that $A=C^{\mathfrak{od}} \otimes_{F} F[\mathbb{Z}/ 2 \mathbb{Z}]$
is a finite dimensional superalgebra with a $\mathbb{Z}/ 2 \mathbb{Z}$-grading induced
by the natural grading of the group algebra $F[\mathbb{Z}/ 2 \mathbb{Z}],$ and
$\mathrm{Id}^{\mathbb{Z}_2}(A)=\mathbb{Z}_2T[\mathrm{Id}(C^{\mathfrak{od}})].$
Thus we have that $\Gamma$ contains
the ideal of graded identities $\mathrm{Id}^{\mathbb{Z}_2}(A)$ of some
finite dimensional $\mathbb{Z}/2 \mathbb{Z}$-graded algebra $A.$
We assume by Lemma \ref{Smu} that $A=\mathcal{O}(A) \times
\mathcal{Y}(A)$ with the senior components $A_1,\dots,A_{\rho}.$
It is clear that $\kappa=\mathrm{ind}_{\mathbb{Z}_2}(\Gamma) \le \kappa_1=\mathrm{ind}_{\mathbb{Z}_2}(A)$
(Lemma \ref{ind1}). If $\Gamma \subseteq \mathrm{Id}^{\mathbb{Z}_2}(A_{i})$ for
some $i=1,\dots,\rho$ then $\kappa_1 =\mathrm{ind}_{\mathbb{Z}_2}(A_i) = \kappa.$
Thus, the case $\Gamma \subseteq
\mathrm{Id}^{\mathbb{Z}_2}(\mathcal{O}(A))$ is trivial.

Assume that $\Gamma \nsubseteq \mathrm{Id}^{\mathbb{Z}_2}(A_i)$ for all
$i=1,\dots,\rho_1$ ($1 \le \rho_1 \le \rho$), and $\Gamma \subseteq
\mathrm{Id}^{\mathbb{Z}_2}(A''),$ where $A''=\times_{i=\rho_1+1}^\rho A_i,$ \
$A'=\times_{i=1}^{\rho_1} A_i.$
Consider the set
$V=S_{\tilde{\mu}}(A') \cap \Gamma$ for
$\tilde{\mu}=\widehat{\mu}(\mathcal{O}(A))$ (Definition \ref{sen}).
Take $\nu=\mathrm{grk}(D)$ for a finitely generated PI-superalgebra
$D$ such that $\Gamma=\mathrm{Id}^{\mathbb{Z}_2}(D).$
By Lemma \ref{GS} there exists a finite
dimensional over $F$ superalgebra $C_\nu$ such
that $\mathrm{Id}^{\mathbb{Z}_2}(C_\nu)=\mathrm{Id}^{\mathbb{Z}_2} \bigl(F\langle X_{\nu}^{\mathbb{Z}_2}
\rangle/\bigl((\mathrm{Id}^{\mathbb{Z}_2}(A')+\mathbb{Z}_2 T[V]) \cap F\langle X_{\nu}^{\mathbb{Z}_2}
\rangle \bigr) \bigr).$

Then $\widetilde{A}=C_\nu \times A'' \times \mathcal{Y}(A)$ is a
finite dimensional superalgebra. For any $f(x_1,\dots,x_n)
\in \mathrm{Id}^{\mathbb{Z}_2}(\widetilde{A}),$ and for any homogenous
polynomials $h_1,\dots,h_n \in F\langle X_\nu^{\mathbb{Z}_2} \rangle$ of the appropriate graded degrees
we have $f(h_1,\dots,h_n)=h+g,$ where $h
\in \mathrm{Id}^{\mathbb{Z}_2}(A'),$ \ $g \in \Gamma \cap
\mathrm{Id}^{\mathbb{Z}_2}(\mathcal{Y}(A)).$ Therefore $h=f(h_1,\dots,h_n)-g \in
\mathrm{Id}^{\mathbb{Z}_2}(A') \cap \mathrm{Id}^{\mathbb{Z}_2}(A'') \cap
\mathrm{Id}^{\mathbb{Z}_2}(\mathcal{Y}(A))=\mathrm{Id}^{\mathbb{Z}_2}(A) \subseteq \Gamma.$
Hence $f(h_1,\dots,h_n)=h+g \in \Gamma$ for any $h_1,\dots,h_n \in
F\langle X_\nu^{\mathbb{Z}_2} \rangle,$ and by Remark \ref{freefg}
$\mathrm{Id}^{\mathbb{Z}_2}(\widetilde{A}) \subseteq \Gamma.$ Particularly
$\mathrm{ind}_{\mathbb{Z}_2}(\widetilde{A}) \ge \kappa.$

Suppose that $A'' \ne 0.$ Then $\kappa_1=$
$\kappa=\mathrm{ind}_{\mathbb{Z}_2}(\widetilde{A}).$ Thus either
$\mathrm{ind}_{\mathbb{Z}_2}(C_\nu)=\kappa$ implying
$\mathcal{O}(\widetilde{A})=A'' \times \mathcal{O}(C_\nu),$ or
$\mathrm{ind}_{\mathbb{Z}_2}(C_\nu) <
\kappa$ implying
$\mathcal{O}(\widetilde{A})=A''.$ Since $S_\mu(A'') \cap \Gamma=\emptyset$ then in the
first case we use Lemmas \ref{dirsumKP}, \ref{Smu}, \ref{subset} and conclude
for any $\mu \ge \max \{ \widehat{\mu}(C_\nu), \tilde{\mu} \}$ we
have $S_\mu(\mathcal{O}(\widetilde{A})) \cap \Gamma=$
$S_\mu(C_\nu) \cap \Gamma \subseteq$ $S_\mu(A') \cap \Gamma
\subseteq$ $V \subseteq \mathrm{Id}^{\mathbb{Z}_2}(C_\nu).$ Thus, in both
cases $S_\mu(\mathcal{O}(\widetilde{A})) \cap \Gamma =
\emptyset,$ and $\widetilde{A}$ satisfies the claims of the
lemma.

If $A'' = 0$ then $\kappa \le \mathrm{ind}_{\mathbb{Z}_2}(\widetilde{A})=$
$\max \{ \mathrm{ind}_{\mathbb{Z}_2}(C_\nu), \mathrm{ind}_{\mathbb{Z}_2}(\mathcal{Y}(A)) \}
\le$ $\kappa_1.$  If $\mathrm{ind}_{\mathbb{Z}_2}(C_\nu) =
\kappa_1=\mathrm{ind}_{\mathbb{Z}_2}(A')$ then $\mathcal{O}(\widetilde{A})=
\mathcal{O}(C_\nu),$ and by analogy with previous case we apply Lemmas \ref{Smu},
\ref{subset} for any $\mu \ge \max \{
\widehat{\mu}(C_\nu), \tilde{\mu} \}$ and obtain
$S_\mu(\mathcal{O}(\widetilde{A})) \cap \Gamma = \emptyset.$
It also follows from Lemma \ref{subset}  that
$\mathrm{ind}_{\mathbb{Z}_2}(\widetilde{A}) = \mathrm{ind}_{\mathbb{Z}_2}(\Gamma),$ and
$\widetilde{A}$ is also the desired algebra.

The last case $A'' = 0,$ $\mathrm{ind}_{\mathbb{Z}_2}(C_\nu) < \kappa_1$ gives
$\widetilde{A}=C_\nu \times \mathcal{Y}(A)$ with
$\mathrm{Id}^{\mathbb{Z}_2}(\widetilde{A}) \subseteq \Gamma,$ and $\kappa \le
\mathrm{ind}_{\mathbb{Z}_2}(\widetilde{A}) < \kappa_1=\mathrm{ind}_{\mathbb{Z}_2}(A).$ Then by
the inductive step on $\mathrm{ind}_{\mathbb{Z}_2}(A)$ we can assume that the
assertion of the lemma holds in this case. \hfill $\Box$

Lemma \ref{Main} with Lemma \ref{subset} implies the
corollary.
\begin{lemma} \label{SmuG}
Let $F$ be a field of characteristic 0.
Let $\Gamma$  be a non-trivial ideal of graded identities of
a finitely generated associative PI-superalgebra over $F.$
Then $\Gamma$ contains the ideal of graded identities
of a finite dimensional associative
$F$-superalgebra $A$ satisfying
$\mathrm{ind}_{\mathbb{Z}_2}(\Gamma)=\mathrm{ind}_{\mathbb{Z}_2}(A),$ and $S_{\tilde{\mu}}(A) =
S_{\tilde{\mu}}(\Gamma)$ for some $\tilde{\mu} \in \mathbb{N}_0$.
\end{lemma}
\noindent {\bf Proof.} Let us take a superalgebra $A$
satisfying the claims of Lemma \ref{Main}.
Lemma \ref{subset} implies that $S_{\mu}(A) \supseteq
S_{\mu}(\Gamma)$ for any $\mu \in \mathbb{N}_0.$ On the
other hand $\mathrm{ind}_{\mathbb{Z}_2}(\Gamma)=\mathrm{ind}_{\mathbb{Z}_2}(A)=(\beta;\gamma).$
Consider an integer $\tilde{\mu}=\max \{ \widehat{\mu}(A), \hat{\mu} \},$ where
$\hat{\mu}$ is defined by Lemma \ref{Main}, and
$\widehat{\mu}(A)$ is taken from Definition \ref{sen}.
Then any polynomial $f \in S_{\tilde{\mu}}(A)$ has the type $(\beta;\gamma-1;\tilde{\mu}),$
and satisfies $f \notin \Gamma$ (since $f \in S_{\tilde{\mu}}(\mathcal{O}(A))$).
Thus, $f \in S_{\tilde{\mu}}(\Gamma).$ \hfill $\Box$

\begin{theorem} \label{GPI}
Let $F$ be a field of characteristic zero.
For any finitely generated associative $\mathbb{Z}/2 \mathbb{Z}$-graded PI-algebra $D$ over $F$
there exists a finite dimensional over $F$ associative superalgebra $\widetilde{C}$ such
that the $\mathbb{Z}_2$T-ideals of $\mathbb{Z}/2 \mathbb{Z}$-graded identities
of $D$ and $\widetilde{C}$ coincide.
\end{theorem}
\noindent {\bf Proof.} Let $\Gamma$ be the ideal of graded identities of $D.$
We use the induction on the Kemer index
$\mathrm{ind}_{\mathbb{Z}_2}(\Gamma)=\kappa=(\beta;\gamma)$ of $\Gamma.$

\underline{Inductive basis.} If
$\mathrm{ind}_{\mathbb{Z}_2}(\Gamma)=(0,0;\gamma)$
then $D$ is a nilpotent finitely generated $F$-superalgebra.
Hence $D$ is finite dimensional.

\underline{Inductive hypothesis.}
Lemma \ref{Main}, and Lemma \ref{SmuG} imply
$\Gamma \supseteq \mathrm{Id}^{\mathbb{Z}_2}(A),$ where
$A=\mathcal{O}(A)+\mathcal{Y}(A)$ is a finite dimensional superalgebra with $\mathrm{ind}_{\mathbb{Z}_2}(\Gamma)=\mathrm{ind}_{\mathbb{Z}_2}(A)=\kappa.$ Moreover,
$S_{\tilde{\mu}}(\Gamma)=S_{\tilde{\mu}}(\mathcal{O}(A))=S_{\tilde{\mu}}(A)
\subseteq \mathrm{Id}^{\mathbb{Z}_2}(\mathcal{Y}(A))$ for some $\tilde{\mu} \in
\mathbb{N}_0.$

Let $A_1,\dots,A_{\rho}$ be the senior components of the algebra $A.$
Denote  $(t_{\bar{0}},t_{\bar{1}})=\beta(\Gamma)=\mathrm{dims}_{\mathbb{Z}_2} A_i,$ \
$t=t_{\bar{0}}+t_{\bar{1}};$ \ $\gamma=\gamma(\Gamma)=\mathrm{nd}(A_i)$ (for all
$i=1,\dots,\rho$).
Let us take for any $i=1,\dots,\rho$
the algebra $\widetilde{A}_i=\mathcal{R}_{q_i,s}(A_i)$
defined by (\ref{FRad}) for the senior component $A_i$
with $q_i=\dim_F A_i,$ \  $s=(t+1)(\gamma+\tilde{\mu}).$
$\widetilde{A}_i$ is a finite
dimensional superalgebra. $\Gamma_i=\mathrm{Id}^{\mathbb{Z}_2}(\widetilde{A}_i)=$
$\mathrm{Id}^{\mathbb{Z}_2}(A_i),$ and $\mathrm{dims}_{\mathbb{Z}_2} \widetilde{A}_i=$
$\mathrm{dims}_{\mathbb{Z}_2} A_i=$ $\beta.$ The Jacobson radical
$J(\widetilde{A}_i)=(X_{q_i}^{\mathbb{Z}_2})/I$ of $\widetilde{A}_i$ is nilpotent of class at
most $s=(t+1)(\gamma+\tilde{\mu}),$ where
$I=\Gamma_i(B_i(X_{q_i}^{\mathbb{Z}_2}))+(X_{q_i}^{\mathbb{Z}_2})^s.$ Here the algebra
$B_i$ can be considered
as the semisimple part of $A_i$ and of $\widetilde{A}_i$
simultaneously (Lemma \ref{Aqs}, Lemma \ref{gamma}).
Particularly, $\widetilde{A}_i$ are superalgebras with elementary decomposition.
By Lemma
\ref{VerbId}, and Lemma \ref{Repr} there exists a finite
dimensional $F$-superalgebra $C$ such that
$\mathrm{Id}^{\mathbb{Z}_2}(C)=\mathrm{Id}^{\mathbb{Z}_2}(\overline{\mathcal{T}}_{\nu}
(\widetilde{A},\Gamma)),$ where $\widetilde{A}=\times_{i=1}^\rho
\widetilde{A}_i,$ $\nu=\mathrm{grk}(D).$

Let us denote $\widetilde{D}_{\nu}=F\langle X_{\nu}^{\mathbb{Z}_2}
\rangle/\bigl((\Gamma + K_{\widetilde{\mu}}(\Gamma)) \bigcap
F\langle X_{\nu}^{\mathbb{Z}_2} \rangle \bigr).$ Lemmas
\ref{ind1}, \ref{addKP} imply that
$\mathrm{ind}_{\mathbb{Z}_2}(\widetilde{D}_{\nu}) \le \mathrm{ind}_{\mathbb{Z}_2}(\Gamma +
K_{\tilde{\mu}}(\Gamma)) < \mathrm{ind}_{\mathbb{Z}_2}(\Gamma)$. By the inductive step
we obtain $\mathrm{Id}^{\mathbb{Z}_2}(\widetilde{D}_{\nu})=\mathrm{Id}^{\mathbb{Z}_2}(\widetilde{U})$
for a finite dimensional over $F$ superalgebra $\widetilde{U}$.
Lemma \ref{VerbId} yields $\Gamma \subseteq
\mathrm{Id}^{\mathbb{Z}_2}(C \times \widetilde{U}).$

Consider a multilinear polynomial $f(x_1,\dots,x_d) \in
\mathrm{Id}^{\mathbb{Z}_2}(C \times \widetilde{U}).$
Let us take any multihomogeneous
polynomials $h_1, \dots, h_d \in F\langle X_{\nu}^{\mathbb{Z}_2} \rangle$
with $\deg_{\mathbb{Z}_2} h_i=\deg_{\mathbb{Z}_2} x_i$ ($i=1,\dots,d$).  We have
$f(h_1,\dots,h_d)=g+h$ for some multihomogeneous graded
polynomials $g \in \Gamma,$ \ $h \in K_{\tilde{\mu}}(\Gamma).$
Then by Lemma \ref{GS1} we obtain $h=f(h_1,\dots,h_d)-g \in \mathcal{S}
\Gamma + \mathrm{SId}^{\mathbb{Z}_2}(\widetilde{A}_i)$ for any
$i=1,\dots,\rho.$ Hence
$\tilde{h}(x_1,\dots,x_n) \in \mathcal{S} \Gamma +
\mathrm{SId}^{\mathbb{Z}_2}(\widetilde{A}_i)$ also holds for the multilinearization
$\tilde{h}$ of $h$. Lemma \ref{Exact4} implies that $\tilde{h}$ is exact for $A_i$ ($\tilde{h} \in K_{\tilde{\mu}}(\Gamma)$) for any $i=1,\dots,\rho$.

Let us fix any $i=1,\dots,\rho.$ Assume that $\{c_1,\dots,c_{q_i} \}$ is a basis of $A_i$ of
homogeneous in the grading elements chosen in (\ref{basisD}),
(\ref{basisU}) (Lemma \ref{Pierce}), and fix the order of these basic elements.
Suppose that
$\bar{a}=(r_1,\dots,r_{\tilde{s}},b_{\tilde{s}+1},\dots,b_n)$ is
any elementary complete evaluation of $\tilde{h}$ by elements of the algebra
$A_i,$ where $r_j \in J(A_i),$ $b_j \in B_i,$
$\tilde{s}=\gamma-1.$
By Lemma \ref{Gammasub} there exist a polynomial
$\hat{h}_{\tilde{\mu}}(\widetilde{Y}_1,\dots,
\widetilde{Y}_{\tilde{s}+\tilde{\mu}},
\widetilde{X},\widetilde{Z}) \in \mathcal{S} \Gamma +
\mathrm{SId}^{\mathbb{Z}_2}(\widetilde{A}_i)$ of the type $(\beta,\gamma-1,\tilde{\mu})$, and
an elementary evaluation $\bar{u}$ in $A_i$
such that $\hat{h}_{\tilde{\mu}}(\bar{u})=\tilde{h}(\bar{a}).$
Moreover, $\hat{h}_{\tilde{\mu}}$ is alternating in any
$\widetilde{Y}_j$ ($j=1,\dots,\tilde{s}+\tilde{\mu}$), and
all variables from $\widetilde{X}
\bigcup \widetilde{Z}$ are exchanged by semisimple elements.
Then we have
\begin{eqnarray} \label{lhdec}
\alpha_1 \hat{h}_{\tilde{\mu}}=\bigl(
\prod_{d=1}^{\tilde{s}+\tilde{\mu}} \mathcal{A}_{\widetilde{Y}_{d, \bar{0}}} \
\mathcal{A}_{\widetilde{Y}_{d, \bar{1}}} \bigr) \hat{h}_{\tilde{\mu}}=
\sum_j  \bigl( \prod_{d=1}^{\tilde{s}+\tilde{\mu}} \mathcal{A}_{\widetilde{Y}_{d, \bar{0}}} \
\mathcal{A}_{\widetilde{Y}_{d, \bar{1}}} \bigr) \bigl(
\mathfrak{\tilde{s}}_j \tilde{g}_j \bigr) (\mathrm{mod } \
\mathrm{SId}^{\mathbb{Z}_2}(\widetilde{A}_i)),
\end{eqnarray}
where $\widetilde{Y}_d=\widetilde{Y}_{d, \bar{0}} \cup
\widetilde{Y}_{d, \bar{1}};$ \  $\alpha_1 \in F,$ $\alpha_1 \ne 0;$ \ $\tilde{g}_j
\in \Gamma,$ \ $\mathfrak{\tilde{s}}_j \in \mathcal{S}.$ Denote $\{\zeta_1,\dots,\zeta_{\hat{n}} \}$
the variables $\widetilde{Y} \bigcup \widetilde{X} \bigcup
\widetilde{Z}$ of $\hat{h}_{\tilde{\mu}}$ (the first variables are
from $\widetilde{Y}=\bigcup_{d=1}^{\tilde{s}+\tilde{\mu}}
\widetilde{Y}_d$).

Given an element $u_l \in A_i$ of the mentioned above evaluation $\bar{u}$ of the polynomial
$\hat{h}_{\tilde{\mu}}$ defined by Lemma \ref{Gammasub}
we take the element $\bar{x}_{\pi(l) \theta_l}=x_{\pi(l) \theta_l}+I$ of the algebra $\widetilde{A}_i.$
Here $x_{\pi(l) \theta_l} \in X_{q_i}^{\mathbb{Z}_2}$ is a graded variable of graded degree
$\theta_l=\deg_{\mathbb{Z}_2} u_l,$ and $\pi(l)$ is the ordinal number of the element
$u_l=c_{\pi(l)}$ in our basis of $A_i,$ \  $1 \le \pi(l) \le q_i$ ($1 \le l \le \hat{n}$).
Consider the following evaluation of
$\hat{h}_{\tilde{\mu}}(\zeta_1,\dots,\zeta_{\hat{n}})$ in the algebra $\widetilde{A}_i$
\begin{eqnarray} \label{Ai}
&& \zeta_l=\varepsilon_{k_1} \bar{x}_{\pi(l) \theta_l}
\varepsilon_{k_2} \in J(\widetilde{A}_i) \quad \mbox{ if } \ \ \zeta_l
\in \widetilde{Y}, \ \  u_l=c_{\pi(l)} \in \varepsilon_{k_1} A_i
\varepsilon_{k_2}; \nonumber \\ && \zeta_l=u_l
\quad \mbox{ if } \ \zeta_l \in \widetilde{X} \bigcup \widetilde{Z}.
\end{eqnarray}

Suppose that in (\ref{lhdec}) the pure trace polynomial $\mathfrak{\tilde{s}}_j$
depends essentially on $\widetilde{Y}.$ Then we get
$\mathfrak{\tilde{s}}_j|_{(\ref{Ai})}=0,$ because the trace of a radical
element is zero (\ref{Atrace}). If $\mathfrak{\tilde{s}}_j$ does not depend
on $\widetilde{Y}$ then $\bigl(
\prod_{d=1}^{\tilde{s}+\tilde{\mu}} \mathcal{A}_{\widetilde{Y}_{d, \bar{0}}} \
\mathcal{A}_{\widetilde{Y}_{d, \bar{1}}} \bigr) \bigl(
\mathfrak{\tilde{s}}_j \tilde{g}_j \bigr) =$
$\mathfrak{\tilde{s}}_j \tilde{\tilde{g}}_j,$ where
$\tilde{\tilde{g}}_j=\bigl(\bigl(
\prod_{d=1}^{\tilde{s}+\tilde{\mu}} \mathcal{A}_{\widetilde{Y}_{d, \bar{0}}} \
\mathcal{A}_{\widetilde{Y}_{d, \bar{1}}} \bigr) \tilde{g}_j \bigr)
\in \Gamma.$
If $\tilde{\tilde{g}}_j|_{(\ref{Ai})} \ne 0$ in $\widetilde{A}_i$ then one
of degree multihomogeneous components of $\tilde{\tilde{g}}_j$ is a $\tilde{\mu}$-boundary polynomial for
$\widetilde{A}_i.$
And it is not a $\tilde{\mu}$-boundary polynomial for $\Gamma,$ because it belongs to $\Gamma.$
This implies that $S_{\tilde{\mu}}(A) \ne S_{\tilde{\mu}}(\Gamma),$
which contradicts to the properties of $A.$ Therefore, we have
$\tilde{\tilde{g}}_j|_{(\ref{Ai})}=0.$ Thus, in any case
$\hat{h}_{\tilde{\mu}}|_{(\ref{Ai})}=0$ holds in the algebra
$\widetilde{A}_i.$ Hence the evaluation
\begin{eqnarray*}
&&\zeta_l=v_l=\varepsilon_{k_1} x_{\pi(l) \theta_l}
\varepsilon_{k_2} \qquad \mbox{ if } \  \  \zeta_l \in \widetilde{Y}, \ \  u_l=c_{\pi(l)}
\in \varepsilon_{k_1} A_i \varepsilon_{k_2};  \\ &&
\zeta_l=v_l=u_l \qquad \mbox{ if } \ \  \zeta_l \in \widetilde{X} \bigcup
\widetilde{Z}
\end{eqnarray*}
of the polynomial $\hat{h}_{\tilde{\mu}}$ in the algebra
$B_i(X_{q_i}^{\mathbb{Z}_2})$ is equal to $\hat{h}_{\tilde{\mu}}(v_1,\dots,v_{\hat{n}}) \in
I=\Gamma_i(B_i(X_{q_i}^{\mathbb{Z}_2}))+(X_{q_i}^{\mathbb{Z}_2})^s.$
Since $|\widetilde{Y}| < s,$ the polynomial $\hat{h}_{\tilde{\mu}}$ is linear in $\widetilde{Y},$
and the variables from $\widetilde{X} \bigcup \widetilde{Z}$ are
replaced by semisimple elements, then we obtain
$\hat{h}_{\tilde{\mu}}(v_1,\dots,v_{\hat{n}}) \in
\Gamma_i(B_i(X_{q_i}^{\mathbb{Z}_2})).$

Consider the map $\varphi:x_{\pi(l) \theta_l} \mapsto
c_{\pi(l)}$ ($l=1,\dots,|\widetilde{Y}|$), and $\varphi(b)=b$ for any $b \in B_i.$
It is clear that $\varphi$ can be extended to a
graded homomorphism $\varphi:B_i(X_{q_i}^{\mathbb{Z}_2}) \rightarrow A_i.$
Then $\varphi(\hat{h}_{\tilde{\mu}}(v_1,\dots,v_{\hat{n}}))=$
$\hat{h}_{\tilde{\mu}}(\varphi(v_1),\dots,\varphi(v_{\hat{n}}))=$
$\hat{h}_{\tilde{\mu}}(\bar{u})=$ $\tilde{h}(\bar{a})
\in $ $\varphi(\Gamma_i(B_i(X_{q_i}^{\mathbb{Z}_2})))=\Gamma_i(A_i)=(0).$

Therefore $\tilde{h}(\bar{a})=0$ holds in $A_i$ for any elementary
complete evaluation $\bar{a} \in A^n$ containing $\gamma-1$
radical elements. Since $\tilde{h}$ is a multihomogeneous
exact polynomial for $A_i,$ and $\gamma=\mathrm{nd}(A_i)$ then $\tilde{h} \in
\mathrm{Id}^{\mathbb{Z}_2}(A_i).$ Hence $h \in \cap_{i=1}^\rho
\mathrm{Id}^{\mathbb{Z}_2}(A_i)=\mathrm{Id}^{\mathbb{Z}_2}(\mathcal{O}(A)),$ and $h \in
\mathrm{Id}^{\mathbb{Z}_2}(\mathcal{O}(A) \times
\mathcal{Y}(A))=\mathrm{Id}^{\mathbb{Z}_2}(A) \subseteq \Gamma.$ Then
we have $f(h_1,\dots,h_d)=g+h \in \Gamma$ for any multihomogeneous
polynomials $h_1, \dots, h_d \in F\langle X_{\nu}^{\mathbb{Z}_2} \rangle$ of
appropriate graded degrees. The application of Remark \ref{freefg} now
implies that $\mathrm{Id}^{\mathbb{Z}_2}(C \times \widetilde{U}) \subseteq \Gamma.$

Therefore, $\Gamma=\mathrm{Id}^{\mathbb{Z}_2}(C \times \widetilde{U}).$ Theorem is proved.
\hfill $\Box$

Observe that to be a PI-algebra is an essential condition for a finitely generated
superalgebra $D$ in Theorem \ref{GPI}, since a finite dimensional superalgebra
is always a PI-algebra.

\section{Specht problem.}

Let $E=\langle e_i, \  i \in \mathbb{N} |\  e_i e_j = -e_j e_i, \
\forall i, j \rangle$ be the Grassmann algebra of infinite rank
with the canonical $\mathbb{Z}/ 2 \mathbb{Z}$-grading $E=E_{\bar{0}} \oplus E_{\bar{1}}$
($E_{\bar{0}}$ and $E_{\bar{1}}$ are the subspaces of $E$ generated by all
monomials in the generators of even and odd lengths respectively).
Consider for a superalgebra $A=A_{\bar{0}} \oplus A_{\bar{1}}$ the Grassmann envelope
$E(A)=A_{\bar{0}} \otimes E_{\bar{0}} \oplus A_{\bar{1}} \otimes E_{\bar{1}}.$
Observe that the algebra $E(A)$ has the natural $\mathbb{Z}/ 2 \mathbb{Z}$-grading
$E(A)_{\theta}=A_{\theta} \otimes E_{\theta},$ \ $\theta \in \mathbb{Z}/2 \mathbb{Z}.$
The next remark is obvious.

\begin{remark} \label{GrIdId}
If two PI-superalgebras $A$ and $B$ has
the same  $\mathbb{Z}/2 \mathbb{Z}$-graded identities
then they have also the same ordinary non-graded polynomial identities.
The T-ideal $\mathrm{Id}(A)$ is the biggest T-ideal that is contained in
$\mathrm{Id}^{\mathbb{Z}_2}(A).$
\end{remark}

\begin{lemma} [Remark 3.7.6, \cite{GZbook}] \label{GrIdId2}
If two PI-superalgebras $A$ and $B$ have
the same  $\mathbb{Z}/2 \mathbb{Z}$-graded identities
then the algebras $E(A),$ and $E(B)$ have the same ordinary non-graded
polynomial identities.
\end{lemma}
\noindent {\bf Proof.} This is a simple consequence of the fact
that the $\mathbb{Z}_2$T-ideals of superalgebras $A$ and $E(A)$
related by some well understood invertible transformation
(see, e.g., \cite{Kem1}, \cite{Kem5}, \cite{GZbook}).
This transformation is completely algorithmic.
Particularly, if two superalgebras are $\mathbb{Z}_2$PI-equivalent then
their Grassmann envelopes are also $\mathbb{Z}_2$PI-equivalent.
By Remark  \ref{GrIdId} the Grassmann envelopes of these superalgebras also
have the same non-graded identities. \hfill $\Box$

\begin{theorem} \label{Z2fg}
The T-ideal of polynomial identities of an associative
PI-algebra over a field of characteristic zero coincides with
the T-ideal of identities of the Grassmann
envelope of some finitely generated PI-superalgebra.
\end{theorem}
\noindent {\bf Proof.} See Theorem 4.8.2 \cite{GZbook}. \hfill $\Box$

Theorem \ref{GPI} along with Theorem \ref{Z2fg}, and Lemma \ref{GrIdId2}
implies the principle Kemer's classification theorem.

\begin{theorem} \label{GPI1}
The T-ideal of polynomial identities of an associative
PI-algebra over a field of characteristic zero coincides with
the T-ideal of identities of the Grassmann
envelope of some finite dimensional superalgebra.
\end{theorem}

Theorem \ref{GPI1} immediately implies the positive solution of the Specht problem.

\begin{theorem} [Theorem 2.4, \cite{Kem1}] \label{GPI2}
The T-ideal of polynomial identities of an associative
PI-algebra over a  field of characteristic zero is finitely
generated as a T-ideal.
\end{theorem}
\noindent {\bf Proof.} Suppose that $\Gamma$ is a T-ideal that is not finitely
based. Then there exists an infinite sequence of multilinear polynomials $\{
f_i(x_1,\dots,x_{n_i}) \}_{i \in \mathbb{N}} \subseteq \Gamma$
satisfying the conditions $\deg f_i < \deg f_j$ for any $i < j$ and
$f_i \notin T[f_1,\dots,f_{i-1}]$ for any $i \in \mathbb{N}.$
Consider the T-ideals $\Gamma_i$
generated by all consequences of the polynomial $f_i$ of degrees
strictly greater then $n_i=\deg f_i$ ($i=1,2,\dots$), and the
T-ideal $\widetilde{\Gamma}=\sum_{i \in \mathbb{N}} \Gamma_i.$
Then we have $f_i \notin \widetilde{\Gamma}$ for any $i \in \mathbb{N}.$
By Theorem \ref{GPI1}
we obtain $\widetilde{\Gamma}=\mathrm{Id}(E(C))$ for a finite dimensional
superalgebra $C.$

Lemma \ref{Pierce} implies $E(C)=E(B) \oplus E(J),$ where $B$
is the semisimple part of $C$, $J$ is the Jacobson radical of
$C.$ Consider a polynomial $f_k$ from our sequence, such that $\deg f_k=n_k >
\mathrm{nd}(C),$  and
consider an evaluation in $E(C)$ of
$f_k$ of the type $x_i =a_i=c_{\delta_i} \otimes g_{\delta_i},$
where $c_{\delta_i} \in
B_{\delta_i} \bigcup J_{\delta_i},$ \
$g_{\delta_i} \in E_{\delta_i},$  \ $\delta_i=\bar{0},\bar{1}$ ($i=1,\dots,n_k$).
If for any $i=1,\dots,n_k$ the elements
$c_{\delta_i}$ are radical then
$f_k(a_1,\dots,a_{n_k})=0$ in $E(C),$ since $n_k >
\mathrm{nd}(C).$ Suppose that $x_{\hat{l}}=
b_{\delta_{\hat{l}}} \otimes
g_{\delta_{\hat{l}}},$ where $b_{\delta_{\hat{l}}} \in
B_{\delta_{\hat{l}}}$ for some $\hat{l}.$ Then
for any element $\tilde{g}_0 \in E_0$ we obtain
$f_k(\dots,a_{\hat{l}},\dots) \tilde{g}_0=$
$f_k(\dots,b_{\delta_{\hat{l}}} \otimes
(g_{\delta_{\hat{l}}} \cdot \tilde{g}_0),\dots)=$
$f_k(\dots,a_{\hat{l}} \cdot (1_B \otimes \tilde{g}_0),\dots)=0,$
since $f_k(x_1,\dots,x_{\hat{l}} \cdot x_0,\dots,x_{n_k}) \in
\widetilde{\Gamma},$ $x_0 \in X,$ \ $1_B$ is the unit of $B.$ It implies
$f_k(a_1,\dots,a_{n_k})=0.$ Therefore $f_k \in \widetilde{\Gamma},$
that contradicts to the choice of $\widetilde{\Gamma}.$
\hfill $\Box$

\bibliographystyle{amsplain}

\end{document}